\documentclass[11pt,twoside]{article} 

\usepackage{amsmath,latexsym,amssymb,amsthm}

\usepackage{hyperref}
\usepackage{titlesec} 
\usepackage{titling} 

\def\titlecolour{\color{DarkSlateBlue}}
\def\seccolour{\color{Crimson}}
\def\headcolour{\color{DarkGrey}}

\pretitle{\begin{center}\LARGE\bfseries\sffamily\titlecolour}
\preauthor{\begin{center}\Large\sffamily\titlecolour}
\postauthor{\par\end{center}}

\titleformat{\section}{\Large\bfseries\sffamily\seccolour}{\thesection}{1em}{}
\titleformat{\subsection}{\large\bfseries\sffamily\seccolour}{\thesubsection}{1em}{}
\titleformat{\paragraph}[runin]{\bfseries\sffamily\seccolour}{}{1em}{}

\usepackage{geometry} 
\usepackage{graphicx} 
\usepackage[margin=1cm,font=small,labelfont=bf]{caption}
\usepackage{subcaption}
\usepackage[svgnames]{xcolor}               
\usepackage{tikz,tikz-cd}
\usepackage{euler}                             
\usepackage{times}

\usepackage{float}
\usepackage{enumerate}

\usepackage{booktabs} 

\usepackage{fancyhdr}
\usepackage[nottoc,notlof,notlot]{tocbibind} 
\usepackage[titles,subfigure]{tocloft} 

\usepackage[T1]{fontenc}

\hypersetup{colorlinks=true,linkcolor=DarkSlateBlue,citecolor=DarkSlateBlue}

\newcounter{todo}

\makeatletter

\newcommand\listtodoname{List of todos}
\newcommand\listoftodos{%
  \section*{\listtodoname}\@starttoc{tod}}
\makeatother

\DeclareMathAlphabet{\mathcal}{OMS}{cmsy}{m}{n}

\numberwithin{equation}{section}
\numberwithin{figure}{section}
\numberwithin{table}{section}

\newtheoremstyle{theoremsf}
{2ex}
{2ex}
{\itshape}
{}
{\sffamily\seccolour}
{}
{1em}
{}

\newtheoremstyle{definitionsf}
{2ex}
{2ex}
{}
{}
{\sffamily\seccolour}
{}
{1em}
{}

\theoremstyle{theoremsf}
\newtheorem{theorem}{Theorem}[section]
\newtheorem{lemma}[theorem]{Lemma}

\theoremstyle{definitionsf}
 \newtheorem{definition}[theorem]{Definition}
 
 \newtheorem{remark}[theorem]{Remark}

\def\gg{\mathfrak{g}}
\def\su{\mathfrak{su}}

\def\hh{\mathfrak{h}}
\def\nn{\mathfrak{n}}
\def\tt{\mathfrak{t}}
\def\zz{\mathfrak{z}}
\def\u{\mathfrak{u}}

\def\xx{\mathbf{x}}

\def\TT{\mathbb{T}}
\def\R{\mathbb{R}}
\def\C{\mathbb{C}}

\def\J{J} 
\def\j{\mathcal{J}}
\def\cO{\mathcal{O}}
\def\CP{\mathbb{CP}}

\def\bb{c}

\def\diag{\mathop\mathrm{diag}\nolimits}
\def\Fix{\mathop\mathrm{Fix}\nolimits}
\def\tr{\mathop\mathrm{tr}\nolimits}

\def\Im{\mathop\mathrm{Im}\nolimits}

\def\e{\mathsf{e}}
\def\i{\mathsf{i}}
\def\Zbar{\overline{Z}}

\def\defn#1{{\bfseries\itshape #1}}

\def\restr#1{\vrule height1.2ex width.4pt
  depth1.4ex\lower0.8ex\hbox{\scriptsize $\,#1$}}

\def\D{({\bf D})}

\newcommand{\WeylChamber}{%
	\path [fill=Pink!50] (0,0) -- (0,10) .. controls (3,9) and  (6,8) .. (8.66,5) -- (0,0);
	\draw [blue] (0,0) -- (0,10);
	\draw [blue] (8.66,5) -- (0,0);
}

\newcommand\blfootnote[1]{%
  \begingroup
  \renewcommand\thefootnote{}\footnote{#1}%
  \addtocounter{footnote}{-1}%
  \endgroup
}

\title{Non-Abelian momentum polytopes for products of $\CP ^2$}
\author{James Montaldi \& Amna Shaddad \\[6pt]\small University of Manchester}
\date{}

\begin{document}
\maketitle
\thispagestyle{empty}

\vbox to 0pt{\hfill\emph{Dedicated to Darryl Holm on the occasion of his 70th birthday}\par\vss}

{\small
\noindent\hrulefill 

\smallskip

\noindent{\large\sffamily\color{DarkSlateBlue} Abstract}

\medskip

\noindent 
This is the first of two companion papers. The joint aim is to study a generalization to higher dimension of the point vortex systems familiar in 2-D.  In this paper we classify the momentum polytopes for the action of the Lie group SU(3) on products of copies of complex projective 4-space.  For 2 copies, the momentum polytope is simply a line segment, which can sit in the positive Weyl chamber in a small number of ways.  For a product of 3 copies there are 8 different types of generic momentum polytope, and numerous transition polytopes, all of which are classified here. The type of polytope depends on the weights of the symplectic form on each copy of projective space.  In the second paper we use techniques of symplectic reduction to study the possible dynamics of interacting generalized point vortices. 

The results can be applied to determine the inequalities satisfied by the eigenvalues of the sum of up to three 3x3 Hermitian matrices with double eigenvalues. 
\medskip

\noindent \emph{MSC 2010}: 53D20 \\[6pt]
\noindent \emph{Keywords}: Momentum map, convex polyhedra, symplectic geometry, eigenvalue estimates

\noindent\hrulefill 
}

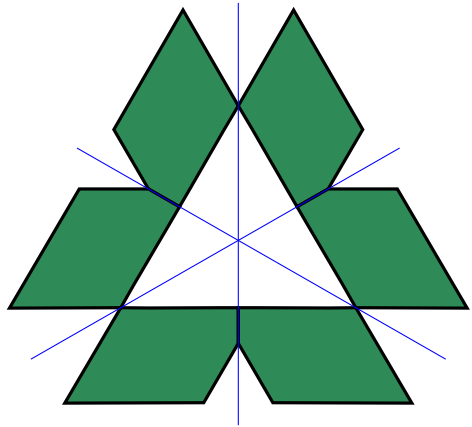
\bigskip
\begin{center}
\begin{tikzpicture}[scale=0.35] 
	\draw[very thick,fill=SeaGreen] (0,5.1) -- (2.078, 8.73) -- (4.676, 4.2)  -- 
	(3.377, 1.95) -- (2.208, 1.276) -- (1.819, 1.95) -- (.5196, 4.22) -- (0,5.1);
	\draw[very thick,fill=SeaGreen] (-0,5.1) -- (-2.078, 8.73) -- (-4.676, 4.2)  -- 
	(-3.377, 1.95) -- (-2.208, 1.276) -- (-1.819, 1.95) -- (-.5196, 4.22) -- (-0,5.1);
	\draw[very thick,fill=SeaGreen,rotate=120] (0,5.1) -- (2.078, 8.73) -- (4.676, 4.2)  -- 
	(3.377, 1.95) -- (2.208, 1.276) -- (1.819, 1.95) -- (.5196, 4.22) -- (0,5.1);
	\draw[very thick,fill=SeaGreen,rotate=120] (-0,5.1) -- (-2.078, 8.73) -- (-4.676, 4.2)  -- 
	(-3.377, 1.95) -- (-2.208, 1.276) -- (-1.819, 1.95) -- (-.5196, 4.22) -- (-0,5.1);
	\draw[very thick,fill=SeaGreen,rotate=-120] (0,5.1) -- (2.078, 8.73) -- (4.676, 4.2)  -- 
	(3.377, 1.95) -- (2.208, 1.276) -- (1.819, 1.95) -- (.5196, 4.22) -- (0,5.1);
	\draw[very thick,fill=SeaGreen,rotate=-120] (-0,5.1) -- (-2.078, 8.73) -- (-4.676, 4.2)  -- 
	(-3.377, 1.95) -- (-2.208, 1.276) -- (-1.819, 1.95) -- (-.5196, 4.22) -- (-0,5.1);
	\draw[thin,blue] (0,-7) -- (0,9);
	\draw[thin,blue,rotate=120] (0,-7) -- (0,9);
	\draw[thin,blue,rotate=-120] (0,-7) -- (0,9);
 \end{tikzpicture}
\end{center}

\newpage

\blfootnote{The picture on the front cover is the intersection of the image of the momentum map with the dual of the Cartan subalgebra, for symplectic weights of type H (see Fig.\,\ref{fig:parameter plane}): it does not appear in the published version. }

\tableofcontents

\section{Introduction}

The now famous convexity theorem of Atiyah, Guillemin and Sternberg and finally Kirwan for the momentum polytope has an interesting history.

In the 1920s, Schur \cite{Schur1923} proved that the diagonal elements $(\delta_1,\ldots,\delta_n)$ of an $n\times n$ Hermitian matrix $A$ satisfy a system of linear inequalities involving the eigenvalues $(\lambda_1,\ldots,\lambda_n)$. In geometric terms, regarding $\delta$ and $\lambda$ as points in $\mathbb{R}^n$ and allowing the symmetric group $S_n$ to act by permutation of coordinates, this result says that $\delta$ lies in the convex hull of the orbit $S_n\cdot\lambda$. 

The converse was proved in the 1950s by Horn \cite{Horn54}, and thus this convex hull is exactly the set of diagonals of the set of all Hermitian matrices with given eigenvalues $(\lambda_1,\ldots,\lambda_n)$.

Kostant generalised these results to any compact Lie group $G$ in the following manner \cite{Kost73}. Consider the coadjoint action of $G$ on the dual $\gg^*$ of its Lie algebra $\gg$. Let $\TT\subseteq G$ be a maximal torus, with Lie algebra $\tt$. Restriction to $\tt$ defines a projection $\gg^*\rightarrow\tt^*$. The Weyl group $W$ acts on $\tt$ and $\tt^*$. Kostant's convexity theorem states,

\noindent\emph{Let $\mathcal{O}\subseteq\gg^*$ be a coadjoint orbit under $G$. Then the projection of $\mathcal{O}$ on $\mathfrak{t}^*$ is the convex hull of a Weyl group orbit.}

The Schur-Horn theorem is the particular case where $G$ is the unitary group $U(n)$ and $\TT$ is the subgroup of diagonal matrices. Then $\gg$ is the Lie algebra of skew-Hermitian matrices. The dual $\gg^*$ can be identified with the set of Hermitian matrices via the pairing  $\left<A,B\right>:=\Im\tr(AB)$, for $A$ Hermitian and $B$ skew-Hermitian. Then the projection of $A\in\gg^*$ on $\tt^*$ is given by the diagonal part of $A$.

This convexity theorem was widely generalised (Atiyah \cite{Atiy82}, Guillemin-Sternberg \cite{GuilSter82}, Kirwan \cite{Kir84b}, etc.). The general relevant framework is that of a symplectic manifold $M$ with a Hamiltonian action of a Lie group $G$. The projection $\mathcal{O}\rightarrow\mathfrak{t}^*$ is a particular case of a momentum map, $M\rightarrow\mathfrak{g}^*$. The most general of these theorems, due to Kirwan, states that the intersection of the image of the momentum map with a positive Weyl chamber in $\tt^*$ is a convex polytope, the \defn{momentum polytope}.

In the vein of the Schur-Horn theorem, this non-Abelian convexity theorem shows for example that if $A$ and $B$ are Hermitian matrices with given eigenvalues, then the eigenvalues of their sum $A+B$ are bounded by a set of linear inequalities involving the given eigenvalues of $A$ and $B$.   See \cite{Knutson} for a description of these ideas.  

In this paper we consider an extended example based on the natural action of $SU(3)$ on $\CP^2$: given $A\in SU(3)$ and $[v]\in\CP^2$ then $A[v]=[Av]$ (in fact the action factors through that of $\mathbb{P}SU(3)$ which is $SU(3)$ factored by the 3-element centre of $SU(3)$, but this has no effect on the material in this paper).  On $\CP^2$ there is an $SU(3)$-invariant symplectic form, the Fubini-Study form, and in fact any invariant symplectic form is a scalar multiple of this particular one.   We consider the compact manifold $M$ given by the product of 2 or 3 copies of $\CP^2$, with the diagonal action of $SU(3)$, and on each copy we choose an invariant symplectic form, with scalars (weights) $\Gamma_j$ (for $j=1,2,3$).  The action of $SU(3)$ on $M$ is Hamiltonian and the momentum map depends on the choice of weights $\Gamma_j$.  The aim of this work is to classify all possible momentum polytopes, depending on the weights.  Note that we use the term `weights' both for the coefficients $\Gamma_j$ and for the weights of a representation: we hope it is clear  from the context which one is meant. In the companion paper \cite{MS19b}, the symplectic weights are called \emph{vortex strengths}. 

The paper is organized as follows.  
After introducing the necessary background in Section\,\ref{sec:Action basic}, Section\,\ref{sec:polytopes2vortices} is dedicated to determining the possible momentum polytopes for the actions of $SU(3)$ on $\CP^2\times\CP^2$, showing there are generically 2 different possible `shapes'; these are just line segments in the positive Weyl chamber (there are also 2 others that are reflections of the first two).   In Section\,\ref{sec:polytopes3vortices} we consider the more interesting case of $\CP^2\times\CP^2\times\CP^2$. We show that, depending on the weights $\Gamma_j$, there are generically 8 distinct types of momentum polytope as well as their reflections under the $*$-involution; there are also numerous transition shapes as the weights vary.  

This is the first of two companion papers; the second \cite{MS19b} uses the results of this paper to study the (reduced) dynamics of a system of generalized point vortices on $\CP^2$, which has symmetry $SU(3)$, acting on a phase space which is the product of copies of $\CP^2$. 
In that paper we discuss the reduced spaces and consider the resulting reduced dynamics and in particular the reduced and relative equilibria and their stability. 

This work forms part of the PhD thesis \cite{AS-thesis}, where further details and alternatives for some of the calculations may be found.

\paragraph{Eigenvalues of Hermitian matrices} 
Following the line of argument of the non-Abelian version of the Schur-Horn theorem mentioned above, one application of our results is to estimating the eigenvalues of the sum of up to three $3\times 3$ Hermitian matrices, each with a double eigenvalue. 

Let $A,B,C$ be three trace-zero $3\times3$ Hermitian matrices each with a double eigenvalue, and let $X=A+B+C$ (if they are not trace zero then replace $A$ by its trace-free part $A_0=A-\frac13\tr(A)\,I_3$, and similarly for $B$ and $C$).   Denote the eigenvalues of $A$ by $\lambda_A,\lambda_A, -2\lambda_A$, and similarly for $B$ and $C$. 

\begin{theorem} 
	\label{thm:eigenvalue estimates}
\begin{enumerate}
\item If\/ $C=0$, then the eigenvalues $\lambda_j$ of\/ $X=A+B$, satisfy $\lambda_1+\lambda_2+\lambda_3=0$ and
$$\lambda_1=\lambda_A+\lambda_B,\quad 
	\lambda_2\in
	\begin{cases}
	[\lambda_A-2\lambda_B,\, \lambda_A+\lambda_B] &\text{ if }\ \lambda_B>0,\\
	[\lambda_A+\lambda_B,\,\lambda_A-2\lambda_B] &\text{ if }\ \lambda_B<0.
	\end{cases}
	$$

\item More generally (with $A,B,C\neq0$), the spectrum of\/ $X=A+B+C$ lies in one of the convex polytopes shown in the figures of Section\,\ref{sec:polytopes3vortices} or its image under the involution $*$, according to the eigenvalues of $A,B,C$.  
\end{enumerate}
Moreover, given any triple $(\lambda_1,\lambda_2,\lambda_3)$ satisfying these inequalities there are Hermitian matrices $A$ with eigenvalues $\lambda_A,\lambda_A,-2\lambda_A$ and $B,C$ with similar eigenvalues, such that $\lambda_1,\lambda_2,\lambda_3$ are the eigenvalues of $A+B+C$. 
\end{theorem}

For example, in part (2), if $\lambda_C=\lambda_B=\lambda_A>0$ then the eigenvalues $\lambda_j$ of $X$ sum to zero and satisfy the inequalities (deduced from Figure\,\ref{subfig:N=3 AAA} and equations\,\eqref{eq:vertices}), 
$$\lambda_j\leq3\lambda_A,\quad j=1,2,3.$$  

Part (1) of this theorem is proved at the end of Section\,\ref{sec:polytopes2vortices}; the proof of part (2) is entirely analogous and is left to the reader.

\section{Hamiltonian action of SU(3) on products of projective spaces}
 \label{sec:Action basic}

In this section we provide the background required, regarding $SU(3)$, symplectic actions and the resulting momentum maps.

\subsection{Background} \label{sec:background}
Recall that if a Lie group acts on a symplectic manifold $(M,\Omega)$ then a momentum map is a map $\J:M\to\gg^*$, where $\gg$ is the Lie algebra and $\gg^*$ its dual vector space, satisfying the differential condition,
\begin{equation}\label{eq:MMdef}
\left<D\J_x(v),\,\xi\right> \ = \ \Omega(\xi_M(x),\,v),
\end{equation}
where $\xi_M$ is the vector field on $M$ associated to $\xi\in\gg$. 

A Lie group $G$ acts naturally on its Lie algebra $\gg$ by the adjoint action and on the dual space $\gg^*$ by the contragredient representation,  the coadjoint action.  An orbit in $\gg^*$ is called a coadjoint orbit.  If, as is our case, the group is compact then the adjoint and coadjoint actions are isomorphic.  If, as we suppose, $G$ is compact and there exists a momentum map, then one can be chosen so that it is equivariant with respect to the given action on $M$ and the coadjoint action on $\gg^*$. 

If $V$ is a symplectic representation of $G$, then the momentum map is given by
$$\left<\J(v),\,\xi\right> \ = \ \tfrac12 [\xi v,\,v],\quad (\xi\in\gg).$$
where $[-,-]$ is the symplectic form. 
An important example is the momentum map for a complex representation $V$ of a torus $\TT$.  Then $V$ is a direct sum of 1-dimensional irreducible representations, of weights $\beta_j,\;j=1,\dots,n$ (where $\dim_\C V=n$, with possible repeats among the $\beta_j$, and possible zeros).  Recall that given a complex representation $V$ of $\TT$, the form $\beta\in\tt^*$ is a weight if the weight-space $V_\beta$ is non-zero, where 
$$V_\beta \ = \ \{v\in V \mid \xi v = i\beta(\xi)v,\;\forall\xi\in\tt\}.$$
If we identify $V$ with $\C^n$, with each coordinate axis being an irreducible representation, then the symplectic form can be written as
$[u,\,v] = \sum_j\Im(u\bar v).$
Then
$$\J(v) \ = \ \tfrac12 \sum_j |v_j|^2\beta_j\ \in \ \tt^*.$$
If instead the symplectic form is altered to $[u,v]=\sum_j\Gamma_j\Im(u_j\overline{v_j})$, then the momentum map becomes
\begin{equation}\label{eq:Gamma_j mom map}
\J(v) \ = \ \tfrac12 \sum_j \Gamma_j|v_j|^2\beta_j\ \in \ \tt^*.
\end{equation}

Coadjoint orbits carry a natural symplectic structure, the Kirillov-Kostant-Souriau, or KKS 2-form, defined by 
$$\omega_{KKS}(\mu)(\xi\cdot\mu,\eta\cdot\mu) := \left<\mu,\,[\xi,\,\eta]\right>,$$
for $\mu\in\gg^*$ and $\xi,\eta\in\gg$. 
If $\cO$ is such a coadjoint orbit with its KKS-form then the coadjoint action of $G$ on $\cO$ is Hamiltonian and the momentum map $\J:\cO\to\gg^*$ is simply given by the inclusion of $\cO$ into $\gg^*$ (see for example \cite{GS-book}).

An important property of momentum maps, often called the bifurcation lemma, and that we will make considerable use of is that, for each $m\in M$,  
\begin{equation}\label{eq:bifurcation lemma}
\mathrm{image}(D\J_m) = \gg_m^\circ,
\end{equation}
where $\gg_m$ is the Lie algebra of the stabilizer of the point $m$, and $\gg_m^\circ$ its annihilator in $\gg^*$. This follows readily from \eqref{eq:MMdef}. 

Given a point $\mu\in\gg^*$ its stabilizer subgroup for the coadjoint action is denoted $G_\mu$.   Any maximal torus $\TT$ of $G_\mu$ is also a maximal torus of $G$.  At the level of Lie algebras, 
$\gg_\mu=\zz_\mu\times \gg'_\mu,$ where $\zz_\mu$ is the centre of $\gg_\mu$. Dualizing, we can write
\begin{equation}\label{eq:z_mu etc}
\gg^* = \zz_\mu^*\times(\gg'_\mu)^*.
\end{equation}
It follows that $\zz_\mu^*=\Fix(G_\mu,\gg^*)$, and similarly we may identify $\gg^*_\mu$ as the subspace of $\gg^*$ given by $\gg_\mu^* := \Fix(\zz_\mu,\gg^*)$. 

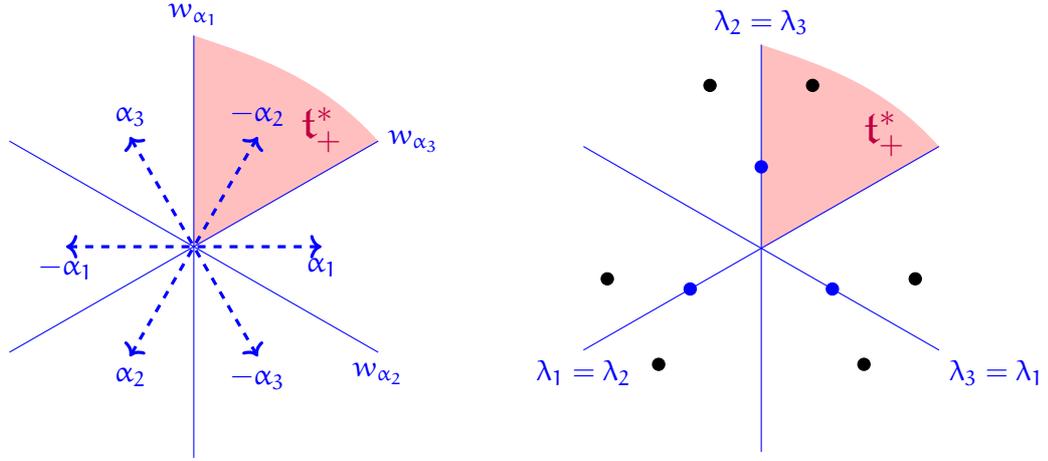
\begin{figure}
    \centering
{ \begin{subfigure} {0.45\textwidth} 
  { \begin{tikzpicture}[scale=0.28]
\path [fill=pink] (0,0) -- (0,10) .. controls (3,9) and  (6,8) .. (8.66,5) -- (0,0);
\draw[purple](6,4) node[anchor=south]{ \LARGE $\mathfrak{t}^*_+$};

\draw [blue] (0,-10) -- (0,10) node[anchor=south]{$w_{\alpha_1}$};
\draw [blue] (-8.66,-5) -- (8.66,5) node [anchor=west]{$w_{\alpha_3}$};
\draw [blue] (-8.66,5) -- (8.66,-5) node[anchor=north]{$w_{\alpha_2}$};

\draw [blue,very thick,dashed,<->] (-6,0) node[anchor=north] {$-\alpha_1$} -- (6,0) node[anchor=north] {$\alpha_1$};
\draw [rotate=-120,blue,very thick,dashed,<->] (-6,0) node[anchor=south] {$-\alpha_2$} -- (6,0) node[anchor=north] {$\alpha_2$};
\draw [rotate=120,blue,very thick,dashed,<->] (-6,0) node[anchor=north] {$-\alpha_3$} -- (6,0) node[anchor=south] {$\alpha_3$};
 \end{tikzpicture}}
 \end{subfigure}}
{ \begin{subfigure}{0.45\textwidth}
  { \begin{tikzpicture}[scale=0.27] \color{black}
\path [fill=pink] (0,0) -- (0,10) .. controls (3,9) and  (6,8) .. (8.66,5) -- (0,0);
\draw[purple](6,4) node[anchor=south]{ \LARGE $\mathfrak{t}^*_+$};

\draw [blue] (0,-10) -- (0,10) node [anchor=south]{$\lambda_2=\lambda_3$};
\draw [blue] (8.66,5) -- (-8.66,-5) node [anchor=north]{$\lambda_1=\lambda_2$};
\draw [blue] (-8.66,5) -- (8.66,-5)node [anchor=north west]{$\lambda_3=\lambda_1$};

\draw [fill] (2.5,8) circle [radius=0.3];
\draw [fill] (-2.5,8) circle [radius=0.3];
\draw [fill] (-7.5,-1.5) circle [radius=0.3];
\draw [fill] (-5,-5.7) circle [radius=0.3];
\draw [fill] (7.5,-1.5) circle [radius=0.3];
\draw [fill] (5,-5.7) circle [radius=0.3];
\draw [fill,blue] (0,4) circle [radius=0.3];
\draw [fill,blue] (3.464,-2) circle [radius=0.3];
\draw [fill,blue] (-3.464,-2) circle [radius=0.3];
 \end{tikzpicture}}
\end{subfigure}}

\caption{On the left the roots for SU(3) and the area shaded in pink is the positive Weyl chamber $\mathfrak{t}^*_+$.  The $\pm\alpha_i$ are the roots.  On the right are shown two orbits of the Weyl group, the black dots show a generic orbit, the blue ones a degenerate orbit.}\label{fig:SU3rootdiagram}
\end{figure}

One particular case is the Cartan subalgebra $\tt$ and its dual $\tt^*=\Fix(\TT,\gg^*)$.  The Weyl group acts on $\tt^*$, and we denote a closed fundamental domain by $\tt^*_+$.  For $SU(3)$, see Figure\,\ref{fig:SU3rootdiagram}.

\paragraph{Conventions} For future computations, especially in Section\,\ref{sec:polytopes3vortices}, we use the following notation and choices.  For a basis of the Cartan subalgebra of $SU(3)$ consisting of diagonal matrices we take,
\begin{equation}\label{eq:basis for su3}
\xi_1=\diag[0,\,i,\,-i],\qquad\xi_2=\diag[-i,0,\,i].
\end{equation}
Here and in what follows, $\diag[a,b,c]$ refers to the matrix with diagonal entries $a,b,c$ and 0s elsewhere.  The positive roots in $\su(3)^*$ are chosen to be $\alpha_1,-\alpha_2$ and $\alpha_3$, where
\begin{equation}\label{eq:roots}
\alpha_1 = \diag[0,1,-1],\quad \alpha_2 = \diag[-1,0,1],\quad \alpha_3 = \diag[1,-1,0]
\end{equation}
(letting $\alpha_2$ be a negative root renders some later expressions more symmetric). See Figure\,\ref{fig:SU3rootdiagram}. 
Notice that we represent elements of $\su(3)$  as skew-hermitian matrices, but elements of the dual $\su(3)^*$ as Hermitian matrices. This requires the pairing to be wrtten as
\begin{equation}\label{eq:pairing}
\mu(\xi)\equiv\left<\mu,\,\xi\right> \ = \ \Im(\tr(\mu\xi)),
\end{equation}
for $\mu\in\su(3)^*$ and $\xi\in\su(3)$. 
Using this, one finds, 
\begin{equation}\label{eq:pairing-values}
\begin{array}{rcl}
\alpha_1(\xi_1) = 2,& \alpha_2(\xi_1)=-1,& \alpha_3(\xi_1)=-1 \\
\alpha_1(\xi_2) = -1,& \alpha_2(\xi_2)=2,& \alpha_3(\xi_2)=-1. 
\end{array}
\end{equation}
With these conventions, the root space decomposition is 
\begin{equation}\label{eq:root spaces}
	\gg=\tt\oplus \gg_{\alpha_1}\oplus \gg_{\alpha_2}\oplus \gg_{\alpha_3},
\end{equation}
where $\gg_{\alpha_1} $ consists of matrices of the form 
$$\xi=\begin{pmatrix}
0&0&0\cr0&0&a\cr0&-\bar{a}&0
\end{pmatrix},
$$
for $a\in\C$, and $ \gg_{\alpha_2}$ (with all entries except $\xi_{13},\xi_{31}$ vanishing) and $ \gg_{\alpha_3}$ similarly.

\paragraph{Witt-Artin decomposition}
Consider a symplectic $G$-manifold $(M,\Omega)$, with $G$-equivariant momentum map
$\J:M\to\gg^*$, and let $m\in M$, and $\mu=\J(m)$.  Then $H=G_m$ acts
symplectically on the tangent space $T_mM$. We recall the Witt-Artin decomposition of $T_mM$, see  \cite{OrteRati04} for details.

Let $T=\gg\cdot m$. It follows from \eqref{eq:MMdef} that
$T^\omega=\ker(D\J_m)$. Consider the four spaces:
\begin{equation}\label{eq:Witt-Artin spaces}
\begin{array}{rclcrcl}
  T_0 &=& T\cap T^\omega\ =\ \gg_\mu.m\,, &\quad& 
  T_1 &=& T/T_0 \ \simeq \ \gg/\gg_\mu\:,\\[8pt]
  N_1 &=& T^\omega/T_0\;, &&
  N_0 &=& T_mM/(T+T^\omega).
\end{array}
\end{equation}
The spaces $T_0$ and $T_1$ give a decomposition of the tangent space $T$ to the group orbit  $G{\cdot}m$ at $m$ while $N_0$ and $N_1$ decompose its (or a) normal
space. $N_1$ is the \defn{symplectic slice} at $m$.

By simple linear algebra, the group action and symplectic form define isomorphisms (of representations of $G_m$),
$$T_0\simeq \gg_\mu/\gg_m,\quad T_1\simeq \gg/\gg_\mu, \quad N_0\simeq T_0^*,$$
and there is a $G_m$-equivariant identification
\begin{equation}\label{eq:Witt-Artin decomposition}
T_mM \simeq T_0 \oplus T_1 \oplus N_1 \oplus N_0.
\end{equation}
In particular, we make a choice for $N_0$ (modulo $T_1$) by requiring $D\J_m(N_0) \subset \gg_\mu^*$, 
which is possible since $D\J_m(T+T^\omega)=D\J_m(T_1)=\gg_\mu^\circ$. With this choice of $N_0$ it follows that 
\begin{equation}\label{eq:N_0 choice}
D\J_m(N_0) = \gg_\mu^*\cap\gg_m^\circ.
\end{equation}
since $\Im(D\J_m)=\gg_m^\circ$. 

If $v\in T_mM$ we write its decomposition with respect to this identification as $v=(w,x,y,z)$, or \begin{equation}\label{eq:wxyz}
v=w + x + y + z \in T_0 \oplus T_1 \oplus N_1 \oplus N_0.
\end{equation}

Finally, $N_1$ and $T_1$ are symplectic while $N_0$ and $T_0$ are isotropic (and paired by the symplectic form). More specifically, given any basis of $T_0$ there is a basis of $N_0$ such that the matrix of $\omega$ at $m$ has the form
$$\left[\omega\right] = \left[\begin{matrix}0&0&0&-I \cr 0&\omega_{T_1}&0&0\cr
0&0&\omega_{N_1}&0\cr I&0&0&0
\end{matrix}\right].$$
Here $\omega_{T_1}$ is the KKS symplectic form on the coadjoint orbit described above, and $\omega_{N_1}$ is the natural symplectic form on the symplectic slice. For details see \cite{OrteRati04}.

\paragraph{MGS normal form} For an action of a compact group $G$ on a manifold $M$, let $m\in M$ and let $S$ be a slice to the orbit (which can be identified with a neighbourhood of 0 in the normal space $N$ to the orbit). Then there is a tubular neighbourhood $U$ of $G\cdot m$ which is equivariantly diffeomorphic to $U\simeq G\times_H S$, where $H=G_m$.

In the symplectic/Hamiltonian setting, this is refined by the Marle-Guillemin-Sternberg normal form, defined as follows, see for example \cite{OrteRati04,Sjamaar} and references therein for details.  The ingredients for this local model are, $\mu\in\tt^*_+$, a closed subgroup $H$ of the stabilizer $G_\mu$ and a symplectic representation $V$ of $H$.  From this one forms a symplectic manifold
$$Y \ = \ Y(\mu,H,V) \ = \ G \times_H (\nn\oplus V),$$
where $\nn=\gg_\mu/\hh$.  The momentum map is given by
\begin{equation}\label{eq:MGS MM}
\J([g,\sigma,v]) = g(\mu+\sigma+\J_{V}(v))g^{-1},
\end{equation}
where $\J_{V}$ is the homogeneous quadratic momentum map for the representation $V$,
$$\left<\J_{V}(v),\,\xi\right> = \tfrac12\,\omega_V(\xi v,v).$$
The momentum polytope for $Y(\mu,H,V)$ is 
$$\Delta(\mu,H,V) = \J(Y)\cap\tt^*_+.$$

Marle and independently Guillemin and Sternberg prove that, given $m\in M$, there is a $G$-invariant neighbourhood $U$ of $m$ and a $G$-invariant neighbourhood $U'$ of $G\times_H(0\times0)$ in $Y(\J(m),G_m,N_1)$ such that $U$ and $U'$ are equivalent as Hamiltonian $G$-spaces. 
Consequently, following Sjamaar \cite{Sjamaar} one makes the following definition:

\begin{definition}\label{def:LMC}
Let $m\in M$ and let $\mu=\J(m)$.  
The \defn{local momentum cone} $\Delta_m$ is defined to be
$$\Delta_m := \Delta(\mu,G_m,N_1)$$
where $N_1$ is the symplectic slice at $m$.  Denote also by $\delta_m$ the germ at $\mu$ of the set $\Delta_m$;  we call $\delta_m$ the \defn{infinitesimal momentum cone} at $m$.
\end{definition}

Sjamaar proceeds to prove the following theorem.

\begin{theorem}[Sjamaar \cite{Sjamaar}]\label{thm:Sjamaar}
Let $M$ be a compact symplectic manifold with a Hamiltonian action of a compact Lie group $G$, and momentum map $\J:M\to\gg^*$. 
\begin{enumerate}
\item If\/ $\J(m_1)=\J(m_2)$ then $\Delta_{m_1}=\Delta_{m_2}$, and \emph{a fortiori} the infinitesimal momentum cones $\delta_{m_1}$ and $\delta_{m_2}$ coincide.
\item The momentum polytope of\/ $M$ is the intersection of all the local momentum cones:
$$\Delta(M) = \bigcap_{m\in\Phi^{-1}(\tt^*_+)}\Delta_m.$$
Moreover, for each $m$, the infinitesimal momentum cone $\delta_m$ coincides with the germ at $\mu$ of $\Delta(M)$.
\item If the point $\mu$ is a vertex of the momentum polytope, then, for any $m\in\J^{-1}(\mu)$,
$$\gg_m^\circ\cap \zz_\mu^*=0.$$
\end{enumerate}
\end{theorem}

Part (3) is stated in a different form by Sjamaar; this equivalent form is proved in \cite{MoRo19}.  In particular, from (3) it follows that a point in the interior of the positive Weyl chamber is a vertex then $G_m=\TT$. 

The statements regarding the infinitesimal momentum cones are not made by Sjamaar, though they are straightforward: by the Marle-Guillemin-Sternberg normal form theorem, there is an invariant neighbourhood of $m$ whose image under the momentum map coincides with a neighbourhood of $\mu$ in $\Delta_m$. Since the momentum map is locally $G$-open onto its image \cite{MT03}, it follows that any sufficiently small representative of the germ $\delta_m$ is a neighbourhood of $\mu$ in $\Delta(M)$.

From (1) we can replace $\Delta_m$ by $\Delta_\mu$ for $\mu=\J(m)$.   It is not always straightforward to find the local momentum cone, although there are 2 cases where it is clear:

\begin{itemize}
\item Firstly, if $m\in\J^{-1}(\mu)$ satisfies $\gg_m=0$ then $D\J_m$ is surjective, and $\Delta_\mu= \tt^*_+$.  
\item Secondly, if $\mu$ lies in the interior of the positive Weyl chamber, then $G_\mu=\TT$ and $G_m$ is a sub-torus.  Then $D\J_m(N_0)=\gg_m^\circ\cap\tt^*$, and moreover, since $G_m$ is a torus, the symplectic representation $N_1$ is a sum of 2-dimensional (symplectic) representations of $G_m$ with  weights $\beta_1,\dots,\beta_r\in\gg_m^*\subset\tt^*$ say. Then (see \eqref{eq:Gamma_j mom map})
$$\J_{N_1}(v_1,\dots,v_r) = \tfrac12\sum_j\Gamma_j|v_j|^2\beta_j \in \gg_m^*,$$
where the coefficients $\Gamma_j$ depend on the symplectic form, 
and it follows that $\Delta_\mu$ is the translation to $\mu$ of the Cartesian product of $\gg_m^\circ\cap\tt^*$ and $\Im(J_{N_1})\subset\gg_m^*$ inside $\tt^*_+$. 
\end{itemize}

\subsection{Momentum map for the SU(3) action on products of \texorpdfstring{$\CP^2$}{CP2}}
\label{sec:MM on CP2s}

We turn our attention to the example of interest, namely $G=SU(3)$ acting on $\CP^2$. Now $\CP^2$ has a particular $SU(3)$-invariant symplectic form known as the Fubini-Study form (obtained from the unit sphere $S^5\subset \C^3$ by reduction by $U(1)$) and denoted $\omega_{FS}$. All other invariant 2-forms on $\CP^2$ are scalar multiples of this basic one. The momentum map for the Fubini-Study form on $\CP^2$ is 
\begin{equation} \label{eq:FS mom-map}
\begin{aligned}
J_0 : \CP^2 &\longrightarrow \su(3)^* \\
Z&\longmapsto Z\otimes\Zbar-\tfrac{1}{3}I_3.
\end{aligned}
\end{equation}
Here $Z=[z_1:z_2:z_3]\in\CP^2$. Since we are viewing $\CP^2$ as the reduction of $S^5$, it follows that $\sum|z_j|^2=1$, and the term $Z\otimes\Zbar$ is the Hermitian matrix $(z_i\bar z_j)$, whose trace is 1 whence the subtraction of the constant term involving the $3\times3$ identity matrix $I_3$. Note that $Z\otimes\Zbar$ is a \emph{Hermitian} matrix, while the elements of $\su(3)$ are skew-Hermitian matrices.   This is not a problem, as the sets of Hermitian and skew-Hermitian matrices are related simply by multiplication by $i$, and re define the pairing of $\su(3)^*$ with $\su(3)$ by the expression in \eqref{eq:pairing}.
It is clear that the expression $J_0$ is equivariant, in that for $g\in SU(3)$,
\begin{equation}\label{eq:equivariance}
J_0(gZ) = g\,J_0(Z)\,g^{-1}.
\end{equation}
In particular, the image of $J_0$ consists of all $3\times 3$ Hermitian matrices with eigenvalues 
$\frac23,\;-\frac13,\;-\frac13.$

The phase space $M$ we are interested in is the Cartesian product of $N$ copies of $\CP^2$, where the $j^{th}$ copy of $\CP^2$ is endowed with an invariant symplectic form $\Gamma_j\omega_{FS}$.  More formally, with $\pi_j:M\to \CP^2$ given by $\pi_j(Z_1,\dots,Z_N) = Z_j$, then the symplectic form on $M$ is 
$$\Omega := \sum_j \Gamma_j\,\pi_j^*\omega_{FS}.$$
We refer to this as the \emph{weighted symplectic form} on $M$, with weights $\Gamma_1,\dots,\Gamma_N$. The momentum map  $\J:M\to\su(3)^*$ for the $SU(3)$-action on $(M,\Omega)$  is then given by 
\begin{equation} \label{eq:MM on CPs}
\J:(Z_1,\ldots,Z_N)\rightarrow\sum^N_{j=1}\Gamma_j\, J_0(Z_J),\quad Z_j\in\CP ^2.
\end{equation}
If is clear from \eqref{eq:equivariance} that this map is also equivariant for the diagonal action on $M$. 

Since it is equivariant, $\J$ descends to a map between orbit spaces we call the \emph{orbit momentum map} and denote $\j$, according to the following diagram,
\begin{equation}\label{eq:OMM}
  \begin{tikzcd}
     M \arrow[r,"\J"]\arrow[d] &\su(3)^*\arrow[d]\\
     M/G\arrow[r,"\j"] &\gg^*/G 
  \end{tikzcd}
\end{equation}
where the vertical maps are the quotient maps. Since every coadjoint orbit in $\gg^*$ intersects $\tt^*$ in a Weyl group orbit, one can identify $g^*/G$ with a positive Weyl chamber $\tt^*_+$.  By the Atiyah-Guillemin-Sternberg-Kirwan convexity theorem, the image $\J(M)/G = \j(M/G)$ is a convex polytope in $\tt^*_+$, called the \emph{momentum polytope}.  For a given number of copies of $\CP^2$, the shape of this polytope will depend on the weights $\Gamma_j$. 

\begin{remark}\label{rmk:quotient v Weyl 1}
While we can identify $\gg^*/G$ with $\tt^*_+$ as described, one needs to be aware that this identification is a homeomorphism but not a diffeomorphism.  Indeed there are many (eg linear) functions on $\tt^*_+$ which are not the restriction of a smooth invariant function on $\gg^*$ (nor Weyl group invariant on $\tt^*$).  See Remark\,\ref{rmk:quotient v Weyl 2} below for further details on this point.
\end{remark}

\begin{remark}\label{rmk:involution}
There is an important involution defined on $\tt^*_+$,  denoted $*$,  defined by
$$*\mu = w(-\mu)$$
where $w$ is the (usually unique) element of the Weyl group that brings $-\mu$ back into the positive Weyl chamber. In the case of $SU(3)$ and the positive Weyl chamber shown in Figure\,\ref{fig:SU3rootdiagram}, we have $w=w_2$. Thus $*\mu=-w_2\mu$; in Figure\,\ref{fig:spectrum diagram} this is the reflection in the line $\lambda_2=0$. The importance in our context is that if one changes $\Gamma=(\Gamma_1,\Gamma_2,\Gamma_3)$ to $-\Gamma$, then 
$$\Delta_\Gamma(M) = *\Delta_{-\Gamma}(M).$$
For this reason it is sufficient to consider $\sum\Gamma_j\geq0$.  
\end{remark}

\subsection{Coadjoint orbits} \label{sec:coadjoint orbits}
We have chosen to represent elements of $\su(3)^*$ as $3\times 3$ Hermitian matrices of trace zero.  The coadjoint action is by conjugation:
$$g\cdot A = gAg^\dagger$$
where $g^\dagger=g^{-1}$ is the conjugate transpose of $g\in SU(3)$. As is well-known from linear algebra courses, two Hermitian matrices are conjugate if and only if they have the same spectrum (including multiplicities). Write this spectrum as $\sigma(A) = \{\lambda_1,\lambda_2,\lambda_3\}$, allowing for multiplicities in the set (sometimes called a multiset). 

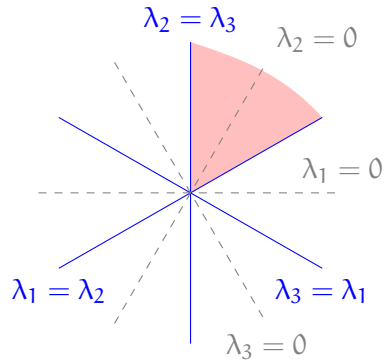
\begin{figure} \centering  %
   \begin{tikzpicture}[scale=0.2] 
\path [fill=pink] (0,0) -- (0,10) .. controls (3,9) and  (6,8) .. (8.66,5) -- (0,0);
\draw [blue] (0,-10) -- (0,10) node [anchor=south]{$\lambda_2=\lambda_3$};
\draw [blue] (8.66,5) -- (-8.66,-5) node [anchor=north]{$\lambda_1=\lambda_2$};
\draw [blue] (-8.66,5) -- (8.66,-5) node [anchor=north]{$\lambda_3=\lambda_1$};
\draw[dashed,gray] (-5, -8.66) -- (5, 8.660) node [anchor=south west] {$\lambda_2=0$};
\draw[dashed,gray] (-10,0) -- (10,0) node [anchor=south] {$\lambda_1=0$};
\draw[dashed,gray] (-5, 8.66) -- (5, -8.660) node [anchor=north] {$\lambda_3=0$};
 \end{tikzpicture}
\caption{This shows the plane parametrized by three real numbers  $\lambda_1,\lambda_2,\lambda_3$ which sum to zero.  The orientation is such that $\lambda_1$ increases to the top of the diagram. Transpositions of the three numbers correspond to reflections in the blue lines. The pink region is where $\lambda_1\geq\lambda_2\geq\lambda_3$. These numbers will be the eigenvalues of a trace zero Hermitian matrix. (Cf.\ the roots shown in Figure\,\ref{fig:SU3rootdiagram})}
\label{fig:spectrum diagram}
\end{figure}

It follows that each coadjoint orbit corresponds to a triple of real numbers summing to zero, and this can be ordered so that $\lambda_1\geq\lambda_2\geq\lambda_3$; see Figure\,\ref{fig:spectrum diagram}. Since the $\lambda_j$ sum to 0, and each is non-zero, it follows that in the preferred ordering, $\lambda_1>0$ and $\lambda_3<0$, while the sign of $\lambda_2$ is variable. In the figure, the coordinate $\lambda_2$ increases as the point moves up or to the left.

\begin{remark}\label{rmk:quotient v Weyl 2}
Continuing Remark\,\ref{rmk:quotient v Weyl 1} above, we note here that the quotient map 
$$\su(3)^*\longrightarrow\su(3)^*/SU(3)$$
can be written as the map $A\mapsto (\chi_2(A),\chi_3(A))$ --- the coefficients of the characteristic polynomial of $A$, which is a smooth and $G$-invariant map. However the map 
$$\su(3)^*\longrightarrow\tt^*_+$$
which maps $A$ to its three eigenvalues, is not smooth but involves extracting roots of the characteristic polynomial. 
\end{remark}

\subsection{Action of $SU(3)$ on products of $\CP^2$}

Let $Z=[z_1:z_2:z_3]\in\CP^2$, then $A\in SU(3)$ acts in a natural way on this point: if $Z'=AZ$ then $z'_j=\sum_k A_{jk}z_k$. Given any $Z\in\CP^2$, the stabilizer $G_Z\simeq U(2)$ is as follows. Consider for example $Z=[1:0:0]$, then $AZ=Z$ if and only if $A$ has the block form
\begin{equation}\label{eq:ISG for one point}
A=\begin{pmatrix}
(\det A_1)^{-1}&0\\
0& A_1
\end{pmatrix}
\end{equation}
where $A_1\in U(2)$.  Since $SU(3)$ acts transitively on $\CP^2$ the stabilizer of any other point will be conjugate to this particular $U(2)$ subgroup. In terms of the root decomposition \eqref{eq:root spaces} of the Lie algebra, this particular copy of $U(2)$ has $\u(2)=\tt\oplus\gg_{\alpha_1}$.  (The stabilizer of $e_j$ has Lie algebra equal to $\tt\oplus\gg_{\alpha_j}$.)

Now consider the diagonal action on $M=\CP^2\times\CP^2$, and let $m=(Z_1,Z_2)\in M$. Let us suppose that $Z_1=e_1=[1:0:0]$.  For $Z_2$ there are 3 cases to consider: first if $Z_1=Z_2$ ( that is, $m$ is a point on the diagonal) then  $G_m$ is again $U(2)$.  Next, if $Z_2$ and $Z_1$ are perpendicular, then we may take $Z_2=e_2=[0:1:0]$ and the stabilizer is the (maximal) torus $\TT^2$ consisting of diagonal matrices:
$$\TT^2 = \left\{\diag[\e^{\i\theta},\,\e^{\i\phi},\,\e^{\i\psi}] \,\mid\, \theta+\phi+\psi=0\mod2\pi\right\}.$$
Finally, if $Z_1,Z_2$ are in general position (neither equal nor perpendicular) then the stabilizer is just a copy of $U(1)$.  For example if $Z_2=[1:1:0]$ then the stabilizer of $(e_1,Z_2)\in M$ is the subgroup of $\TT^2$ consisting of matrices of the form 
$$\left\{\mathrm{diag}\left[\e^{\i\theta},\,\e^{\i\theta},\,\e^{-2\i\theta}\right]\right\} \ \simeq \ U(1).$$  
We summarize these possibilities in the following table,
\begin{equation}\label{eq:table of ISGs N=2}
\begin{tabular}{c|c}
geometry of $m$ & stabilizer \\
\hline
on diagonal & $U(2)$ \\
$(u,u^\perp)$ & $\TT^2$\\
general position & $U(1)$ 
\end{tabular}
\end{equation}
where $u$ and $u^\perp$ are any pair of orthogonal points in $\CP^2$.

For a product of three copies of $\CP^2$ the analysis is similar.  We have:
\begin{equation}\label{eq:table of ISGs N=3}
\begin{tabular}{c|c}
geometry & stabilizer \\
\hline
on diagonal & $U(2)$ \\
$(u,u,v)$ & $\TT^2$\\
$(u,v,w)$ & $\TT^2$ \\
$(u,u',v)$ & $U(1)$ \\
spanning a plane & $U(1)$ \\
general position & $\mathbf{1}$
\end{tabular}
\end{equation}
Here $u, v$ and $w$ are pairwise orthogonal, while $u,u'$ are distinct but not orthogonal, and $v$ is orthogonal to both $u$ and $u'$.  

The following lemma will be useful when computing local momentum cones in Section\,\ref{sec:polytopes3vortices}.  Recall that $e_1=[1:0:0]\in\CP^2$ etc.. Recall also that given a complex representation $V$ of $\TT$, the form $\alpha\in\tt^*$ is a weight if the weight-space $V_\alpha$ is non-zero, where 
$$V_\alpha \ = \ \{v\in V \mid \xi v = i\alpha(\xi)v,\;\forall\xi\in\tt\}.$$

\begin{lemma}\label{lemma:weights}
The representation of\/ $\TT^2$ on the tangent space $T_{e_i}\CP^2$, has the following weights:
$$T_{e_1}\CP^2 = -\alpha_3\oplus \alpha_2,\quad T_{e_2}\CP^2 =  -\alpha_1\oplus\alpha_3,\quad T_{e_3}\CP^2 = -\alpha_2 \oplus\alpha_1.$$
\end{lemma}

See \eqref{eq:roots} for the definition of the $\alpha_j$; the choice of $\pm$-signs in each case is determined by the natural complex structure on $T_{e_j}\CP^2$.  This sign is compatible with the symplectic structure if the corresponding symplectic weight satisfies $\Gamma>0$.  

\begin{proof}
For $T_{e_1}\CP^2$, the tangent vectors are of the form $\xx=(0,v,w)^T$ (with $v,w\in \C$), and the action of $\xi_1$ on this is $\xi_1.(0,v,w)^T=(0,iv,-iw)^T$, and $\xi_2.(0,v,w)^T = (0,iv,2iw)^T$. Note that the action of $\xi_2$ has to be adjusted to ensure it fixes $e_1$; the matrix $\diag[0,i,2i]$ acts the same as the one given in \eqref{eq:basis for su3}, and this one manifestly fixes $e_1$ and hence acts on $T_{e_1}\CP^2$ by simple multiplication.  Thus, $(0,v,0)^T$ has weight satisfying $\alpha(\xi_1)=\alpha(\xi_2)=1$, hence its weight is $-\alpha_3$, while $(0,0,w)^T$ has weight satisfying $\alpha(\xi_1)=-1$ and $\alpha(\xi_2)=2$, giving the weight $\alpha_2$. The other cases are similar. In particular,  $\xx\in T_{e_2}\CP^2$ can be written $\xx=(u,0,w)^T\in\C^2$ and the weight of $(0,0,w)^T$ is $-\alpha_1$ while that of $(u,0,0)^T$ is $\alpha_3$. 
Finally, $\xx\in T_{e_3}\CP^2$ can be written $\xx=(u,v,0)^T\in\C^2$ and the weight of $(u,0,0)^T$ is $-\alpha_2$ while that of $(0,v,0)^T$ is $\alpha_1$. 
\end{proof}

\section{Momentum polytopes for SU(3) action on \texorpdfstring{$\CP ^2\times\CP ^2$}{CP2 x CP2}}\label{sec:polytopes2vortices}

To determine these polytopes one can apply a far simpler argument than for the product of 3 copies, as we shall see.  This example has been considered before by Bedulli and Gori \cite{BG07}. 

The action of $SU(3)$ on $M=\CP^2\times\CP^2$ is not transitive, and it is not hard to see that $(Z_1,Z_2)$ and $(Z_1',Z_2')$ lie in the same orbit if and only if the distance between $Z_1$ and $Z_2$ is equal to that between $Z_1'$ and $Z_2'$.  
It follows that the orbit space $M/SU(3)$ is a compact line segment, parametrized by this distance.  The image of the orbit momentum map $\j$ is therefore 1-dimensional, and by the convexity theorem it must be a line segment (or a point). A line segment has two ends, and it suffices to find these two end points, which will necessarily be the images of the end-points of $M/SU(3)$.

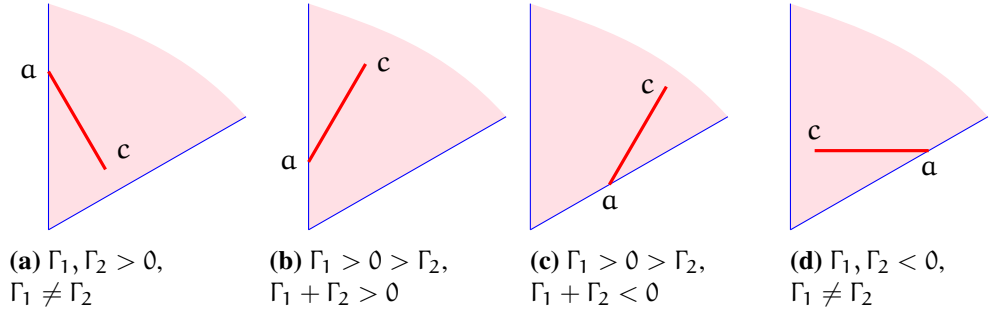
\begin{figure}[t]  
\centering
{\begin{subfigure}{0.2\textwidth}
\begin{tikzpicture}[scale=0.3] 
\WeylChamber
\draw[very thick,red] (0,7) node[anchor=east,black]{$a$} -- (2.500, 2.670) node[black,anchor=south west]{$\bb$};
 \end{tikzpicture}
\caption{$\Gamma_1,\Gamma_2>0$,\\ $\Gamma_1\neq\Gamma_2$}
\label{subfig:generic N=2 a}
\end{subfigure}}
\quad
{\begin{subfigure}{0.2\textwidth}
\begin{tikzpicture}[scale=0.3] 
\WeylChamber
\draw[very thick,red] (0,3) node[anchor=east,black]{$a$} -- (2.500, 7.330) node[black,anchor=west]{$\bb$};
 \end{tikzpicture}
\caption{$\Gamma_1>0>\Gamma_2$, \\ $\Gamma_1+\Gamma_2>0$}
\label{subfig:generic N=2 b}
\end{subfigure}}
\quad
{\begin{subfigure}{0.2\textwidth}
\begin{tikzpicture}[scale=0.3] 
\WeylChamber
\draw[very thick,red] (3.464, 2.) node[anchor=north,black]{$a$}-- (5.964, 6.330) node[black,anchor=east]{$\bb$};
 \end{tikzpicture}
\caption{$\Gamma_1>0>\Gamma_2$, \\ $\Gamma_1+\Gamma_2<0$}
\label{subfig:generic N=2 c}
\end{subfigure}}
\quad
{\begin{subfigure}{0.2\textwidth}
\begin{tikzpicture}[scale=0.3] 
\WeylChamber
\draw[very thick,red] (6.062, 3.5) node[anchor=north,black]{$a$} -- (1.062, 3.500) node[black,anchor=south]{$\bb$};
 \end{tikzpicture}
\caption{$\Gamma_1,\Gamma_2<0$, \\$\Gamma_1\neq\Gamma_2$}
\label{subfig:generic N=2 d}
\end{subfigure}}
\bigskip

\caption{The four generic polytopes for the action of $SU(3)$ on $\CP^2\times\CP^2$. In each case $a$ represents the image of  points of the form $(u,u)$, and $\bb$ of points of the form $(u,u^\perp)$. Notice that all these poytope-segments are parallel to one of the roots (equivalently, orthogonal to one of the walls of the Weyl chamber). Notice that figures (a) and (d) are related by the involution $*$ of Remark\,\ref{rmk:involution}, as are figures (b) and (c).}
\label{fig:generic N=2}
\end{figure}

\begin{figure}  
\centering
\begin{subfigure}{0.2\textwidth}
\setcounter{subfigure}{4}
\begin{tikzpicture}[scale=0.3] 
\WeylChamber
\draw[very thick,red] (0,7) node[anchor=east,black]{$a$} -- (3.031, 1.750) node[black,anchor=north]{$\bb$};
 \end{tikzpicture}
 \caption{$\Gamma_1=\Gamma_2>0$}
   \label{subfig:transitional N=2 a}
\end{subfigure}
\qquad
{\begin{subfigure}{0.2\textwidth}
\begin{tikzpicture}[scale=0.3] 
\WeylChamber
\draw[very thick,red] (0,0) node[anchor=east,black]{$a$} -- (2.500, 4.330) node[black,anchor=west]{$\bb$};
 \end{tikzpicture}
 \caption{$\Gamma_1+\Gamma_2=0$}
  \label{subfig:transitional N=2 b}
 \end{subfigure}}
\qquad
{\begin{subfigure}{0.2\textwidth}
\begin{tikzpicture}[scale=0.3] 
\WeylChamber
\draw[very thick,red] (6.062, 3.5) node[anchor=north,black]{$a$} -- (0., 3.500) node[black,anchor=east]{$\bb$};
 \end{tikzpicture}
 \caption{$\Gamma_1=\Gamma_2<0$}
  \label{subfig:transitional N=2 c}
\end{subfigure}}
\bigskip

\caption{The three transitional polytopes for the action of $SU(3)$ on $\CP^2\times\CP^2$. See the caption of Figure\,\ref{fig:generic N=2} for explanations of notation, and Remark\,\ref{rmk:transitions 2} for discussion. Note that the involution $*$ exchanges figures (e) and (g) and leaves (f) unchanged.}
\label{fig:transitional N=2}
\end{figure}
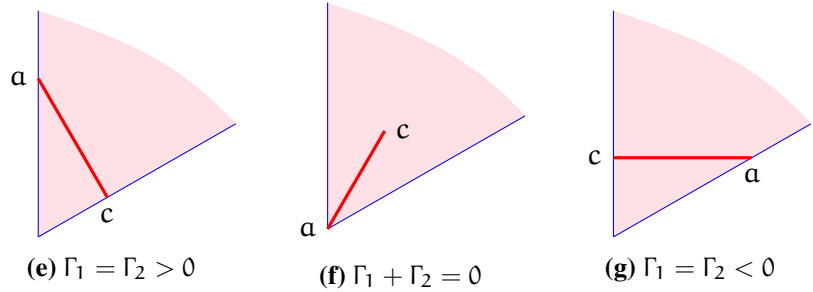

\begin{theorem}\label{thm:2-polytopes} 
The momentum polytopes $\Delta(M)$ of the $SU(3)$ action on $M=\CP ^2\times\CP ^2$ with weighted symplectic form $\Gamma_1\omega_{FS}\oplus\Gamma_2\omega_{FS}$ with $\Gamma_i\neq0$ fall into four different types for which $\Gamma_1\neq\pm\Gamma_2$, and three transitional ones where $\Gamma_1=\pm\Gamma_2$; these are shown in Figures\,\ref{fig:generic N=2} and \ref{fig:transitional N=2} respectively.
\end{theorem} 

\begin{remark}\label{rmk:types of polytope}
We have not defined what we mean by the \emph{type} of a momentum polytope. Without giving a formal definition, the type is a combination of `combinatorics' and `geometry' within the positive Weyl chamber. For example in Figure\,\ref{fig:generic N=2}, (a) and (b) have the same `combinatorics' (indeed all 4 figures do), but their geometry relative to the Weyl chamber is different.  
\end{remark}

\begin{proof}
First let us find the images of the momentum map for the points in $M$ at extrema of the distance.  It suffices to choose representatives: for minimal distance we have $m=(e_1,e_1)$ with $e_1=[1:0:0]$ and for maximal distance we chose points that are perpendicular, for example $(e_1,e_2)$, with $e_2=[0:1:0]$.  From \eqref{eq:FS mom-map} and \eqref{eq:MM on CPs} one finds,
$$\J(e_1,e_2) = \tfrac{1}3\Gamma_1
\begin{pmatrix}
2 &0&0\cr 0& -1&0\cr 0&0&-1
\end{pmatrix}
+ \tfrac{1}3\Gamma_2\begin{pmatrix}
-1 &0&0\cr 0& 2&0\cr 0&0&-1
\end{pmatrix}.
$$
From this and a similar calculation for $\J(e_1,e_1)$ we obtain the spectra,
\begin{eqnarray*}
\sigma(\J(e_1,e_1)) &=& \left\{\tfrac{2(\Gamma_1+\Gamma_2)}3,\,\tfrac{-(\Gamma_1+\Gamma_2)}3,\, \tfrac{-(\Gamma_1+\Gamma_2)}3\right\}, \\ 
\sigma(\J(e_1,e_2)) &=& \left\{\tfrac{(2\Gamma_1-\Gamma_2)}3,\,\tfrac{(2\Gamma_2-\Gamma_1)}3,\, \tfrac{-(\Gamma_1+\Gamma_2)}3\right\}.
\end{eqnarray*}

Notice that the first of these spectra has two equal eigenvalues so lies on a line of reflection. When ordered by decreasing value, this point is marked $a$ in each diagram. All three are equal if and only if $\Gamma_1+\Gamma_2=0$, and in that case $a$ lies at the origin.  The other point $\bb=\J(e_1,e_2)$  generically has 3 distinct elements, so does not lie on a line of reflection.  Repeated eigenvalues occur if and only if $\Gamma_1=\Gamma_2$, as is readily checked (or if one of the $\Gamma_j$ vanishes, which we are excluding), and the corresponding point is marked $\bb$ in the two figures.

There remains to show that the points $a,\bb$ are indeed the endpoints of the segment $\Delta(M)$ as claimed.  Clearly, since $a$ lies on a wall of the Weyl chamber, it must be an endpoint of the segment.   Now any endpoint is a vertex of the polytope so corresponds either to a point in the wall or a fixed point for the torus action, but $a$ and $\bb$ are the only images of fixed points.
\end{proof}

\begin{remark}\label{rmk:transitions 2}
As $(\Gamma_1,\Gamma_2)$ varies in the plane there are transitions that occur at the points described in the theorem. Here we briefly describe these.

First, as $\Gamma_1+\Gamma_2$ goes from being positive to negative, the transition from Figure\,\ref{subfig:generic N=2 b} to \ref{subfig:generic N=2 c} is seen clearly through Figure\,\ref{subfig:transitional N=2 b}. When $\Gamma_1+\Gamma_2=0$ one of the eigenvalues vanishes for both points $a$ and $\bb$ so the segment lies along the line $\lambda_2=0$. 

The transition between Figures\,\ref{subfig:generic N=2 a} and \ref{subfig:generic N=2 b} occur as $\Gamma_2$ changes sign.  As $\Gamma_2\to0$, the segment in Figure\,\ref{subfig:generic N=2 a} or \ref{subfig:generic N=2 b} becomes shorter, and in the limit becomes just the point $a$ (when the symplectic form is degenerate, the momentum map does not `see' the second factor in the product $M$, and the momentum polytope reduces to that of $\CP^2$ which is just a single point). The transition between Figures\,\ref{subfig:generic N=2 c} and \ref{subfig:generic N=2 d} is similar.

Finally, the transitional figures shown in Figures\,\ref{subfig:transitional N=2 a} and \ref{subfig:transitional N=2 c} occur when $\Gamma_1=\Gamma_2$ and $\J(e_1,e_2)$ has a double eigenvalue. As say, $\Gamma_1$ decreases through the value $\Gamma_2$ from $\Gamma_1>\Gamma_2>0$ to $\Gamma_2>\Gamma_1>0$, the segment in Figure\,\ref{subfig:generic N=2 a} extends until it hits the right-hand wall (as in Fig.\,\ref{subfig:transitional N=2 a}) and then retreats back to look like the segment in Figure\,\ref{subfig:generic N=2 a} again.
\end{remark}

\begin{proof}[Proof of Theorem\,\ref{thm:eigenvalue estimates} (Part 1)]
Consider $M=\CP^2\times\CP^2$, with $\Gamma_1=-3\lambda_A$ and $\Gamma_2=-3\lambda_B$. 
The two extremes of the segment have spectra given above 
\begin{eqnarray*}
\sigma(\J(e_1,e_1)) &=& \left\{\lambda_A+\lambda_B,\; \lambda_A+\lambda_B,\;-2(\lambda_A+\lambda_B)\right\}, \\ 
\sigma(\J(e_1,e_2)) &=& \left\{\lambda_A+\lambda_B,\;\lambda_A-2\lambda_B,\;\lambda_B-2\lambda_A\right\}.
\end{eqnarray*}
Then $\lambda_1=\lambda_A+\lambda_B$, while $\lambda_2$ lies between $\lambda_A+\lambda_B$ and $\lambda_A-2\lambda_B$ as stated in the theorem.
\end{proof}

\section{Momentum polytopes for SU(3) action on \texorpdfstring{$\CP ^2\times\CP ^2\times\CP ^2$}{CP2 x CP2 x CP2}}
\label{sec:polytopes3vortices}

Now let $M=\CP^2\times\CP^2\times\CP^2$. 
Recall from Section\,\ref{sec:MM on CP2s}, the momentum map for the $SU(3)$ action on the manifold $M$ is
\begin{equation*}
\J:(Z_1,Z_2,Z_3)\longmapsto \sum_{j=1}^3\Gamma_j \, Z_j\otimes \Zbar_j - \tfrac{1}{3}\left(\sum_{j=1}^3\Gamma_j\right)I_3
\end{equation*}
with $Z_1,Z_2,Z_3\in\CP ^2$, see \eqref{eq:FS mom-map}.

\begin{figure}[p]
	\centering


  \begin{tikzpicture}[scale=1.3]   
		\def\blueline{[blue] (-4.5,0) -- (4.5,0)}
	\draw [yshift=2cm] \blueline;
	\draw [yshift=5mm] \blueline;
	\draw [yshift=-1cm,dashed] \blueline;
	\draw [rotate=120,yshift=2cm] \blueline;
	\draw [rotate=120,yshift=5mm] \blueline;
	\draw [rotate=120,yshift=-1cm,dashed] \blueline;
	\draw [rotate=-120,yshift=2cm] \blueline;
	\draw [rotate=-120,yshift=5mm] \blueline;
	\draw [rotate=-120,yshift=-1cm,dashed] \blueline;
	\draw [rotate=90,scale=1.1] \blueline;
	\draw [rotate=210,scale=1.1] \blueline;
	\draw [rotate=-30,scale=1.1] \blueline;
	\draw[ultra thick,dashed] (0,2) -- (1.732,-1) -- (-1.732,-1) -- (0,2);

	\draw[line width=2.5pt,red] (0,0) -- (4.5,2.598);
	\draw[line width=2.5pt,red] (0,0) -- (4.5,-2.598);

	\draw [fill] (0.3,0) circle [radius=0.05] node [anchor=west]{$A$};
	\draw [fill] (0.8,-0.2) circle [radius=0.05] node [anchor=west]{$B$};
	\draw [fill] (1.7,-0.2) circle [radius=0.05] node [anchor=west]{$C$};
	\draw [fill] (3.5,-1.6) circle [radius=0.05] node [anchor=west]{$D$};
	\draw [fill] (3.5,-0.2) circle [radius=0.05] node [anchor=west]{$E$};
	\draw [fill] (3.5,1.2) circle [radius=0.05] node [anchor=west]{$F$};
	\draw [fill] (4.3,2.2) circle [radius=0.05] node [anchor=west]{$G$};
	\draw [fill] (2,0.8) circle [radius=0.05] node [anchor=west]{$H$};

	\draw[yshift=4pt,blue] (4.5,-1) node {\small$\Gamma_2=0$};
	\draw[rotate=120,yshift=-4pt,blue] (4.5,-1) node[rotate=-60] {\small$\Gamma_3=0$};
	\draw[rotate=240,yshift=4pt,blue] (4.5,-1) node[rotate=60] {\small$\Gamma_1=0$};
	\draw[yshift=10pt,red] (2.2,1.2) node[rotate=30] {\small$\Gamma_1=\Gamma_2$};
	\draw[yshift=-10pt,red] (3,-1.6) node[rotate=-30] {\small$\Gamma_2=\Gamma_3$};
	\draw[blue] (0,4) node [anchor=north,rotate=90] {\small$\Gamma_1=\Gamma_3$};
	\draw[yshift=4pt,blue] (1,2) node {\small$\Gamma_1+\Gamma_3=0$};
	\draw[rotate=-120,yshift=4pt,blue] (1,2) node [rotate=60]{\small$\Gamma_2+\Gamma_3=0$};
	\draw[rotate=120,yshift=4pt,blue] (1,2) node [rotate=-60]{\small$\Gamma_1+\Gamma_2=0$};
	\draw[yshift=4pt,blue] (4.5,0.5) node {\small$\Gamma_2=\Gamma_1+\Gamma_3$};
	\draw[rotate=240,yshift=6pt,blue] (4.5,0.5) node[rotate=60] {\small$\Gamma_1=\Gamma_2+\Gamma_3$};
	\draw[rotate=120,yshift=6pt,blue] (4.5,0.5) node[rotate=-60] {\small$\Gamma_3=\Gamma_1+\Gamma_2$};
  \end{tikzpicture}
\caption{This shows the parameter plane $\Gamma_1+\Gamma_2+\Gamma_3=\text{const}$ with const${}>0$.  Within the central black triangle all 3 weights are positive. The value of $\Gamma_2$ is constant on horizontal lines and increases vertically upwards; variations of the other variables can be deduced from this.  The blue lines indicate where the polytope type changes, see Table\,\ref{table:transitions}. The sector between the red lines is where $\Gamma_1\geq\Gamma_2\geq\Gamma_3$.   The generic polytope types are labelled $A,B,\dots,H$, and illustrated in Fig.\,\ref{fig:polytopes A-H}, and the respective transitions are labelled AB, CE etc., see Fig.\,\ref{fig:transition parameters}.}
\label{fig:parameter plane}
\end{figure}
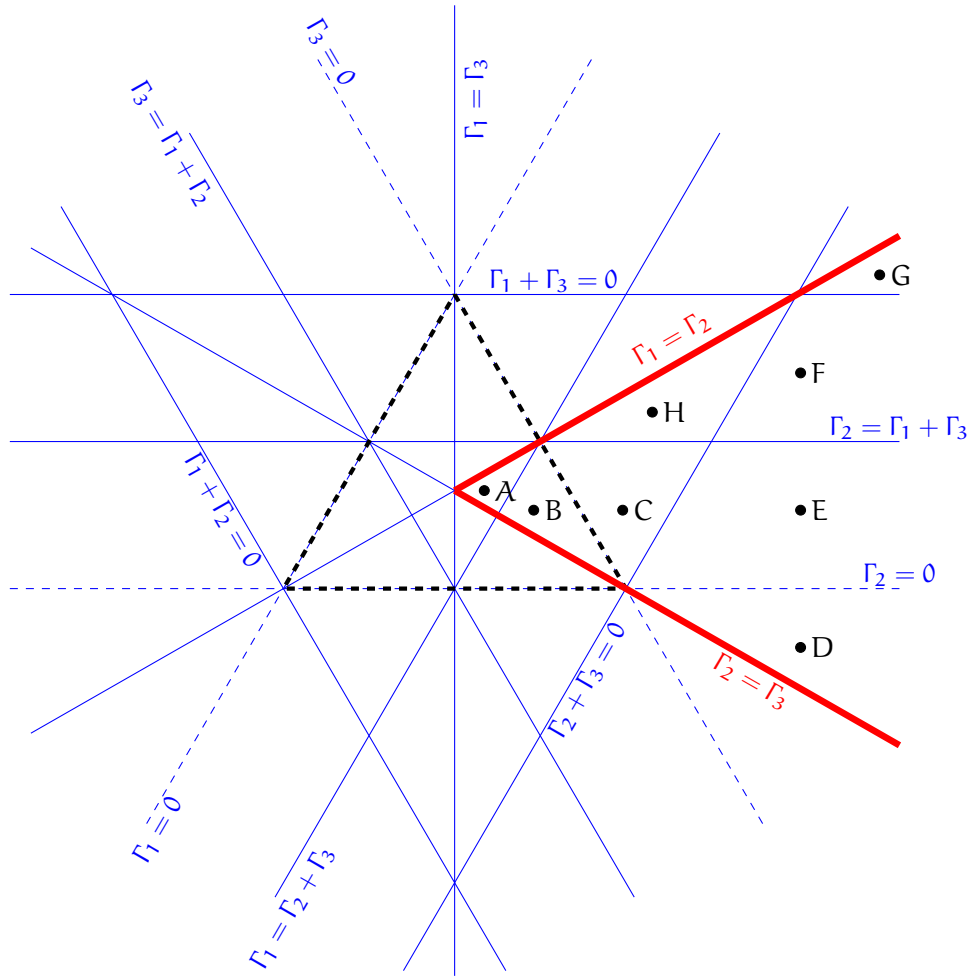

\begin{figure}[p] 
\centering

\centering

{\begin{subfigure}{0.35\textwidth}
\begin{tikzpicture}[scale=0.4] 
\WeylChamber
	\draw[very thick,fill=Crimson] (0, 8.25) -- (3.572, 2.062) -- 
	 (2.598, 1.5) -- 	(1.299, 1.5) -- (0, 3.75) -- (0, 8.25);
	\draw[thick] (1.299, 6) -- (0, 3.75);
	\draw[thick] (2.598, 3.75) -- (1.3,1.5);
	\draw[thick] (3.248, 2.625) -- (2.598, 1.5);
	 \draw [fill] (0, 8.25) node[anchor=east]{$a$};
	 \draw [fill] (1.3, 1.5) node[anchor=north]{$b$};
	 \draw [fill] (1.299, 6) circle [fill,radius=0.08] node[anchor=south west]{$c_3$};
	 \draw [fill] (2.598, 3.75) circle [fill,radius=0.08] node[anchor=south west]{$c_2$};
	 \draw [fill] (3.248, 2.625) circle [fill,radius=0.08] node[anchor=south west]{$c_1$};
 \end{tikzpicture}
\caption{Polytope A}
\label{subfig:N=3 A}
\end{subfigure}}
\qquad
{\begin{subfigure}{0.35\textwidth}
\begin{tikzpicture}[scale=0.4] 
\WeylChamber
	\draw[very thick,fill=Crimson] (0,8.25) node[anchor=east]{$a$} -- 
	(1.3,6) circle [fill,radius=0.08] node[black,anchor=south west]{$c_3$} --
	(1.948, 4.875) circle [fill,radius=0.08]  node[black,anchor=south west]{$c_2$} --
	(3.248, 2.625) node[black,anchor=west]{$c_1$} --
	(0.6495, 2.625) node[black,anchor=north]{$b$} --
	(0,3.75) -- (0,8.25);
	\draw [thick] (1.948, 4.875) -- (0.6495, 2.625);
	\draw [thick] (1.3,6) -- (0,3.75);
 \end{tikzpicture}
\caption{Polytope B}
\label{subfig:N=3 B}
\end{subfigure}}

\vskip 5mm

{\begin{subfigure}{0.35\textwidth}
\begin{tikzpicture}[scale=0.4] 
\WeylChamber
	\draw[very thick,fill=Crimson] (0, 6.75) node[anchor=east]{$a$} -- 
	(1.3, 9) node[black,anchor=west]{$c_3$} --
	(3.248, 5.625) node[black,anchor=west]{$b$} --
	(1.948, 3.375) node[black,anchor=north]{$c_2$} --
	(0.6495, 5.625) circle[radius=0.08] node[black,anchor=north]{$c_1\;$} -- (0,6.75);
	\draw[thick] (3.248, 5.625) -- (0.6495, 5.625);
 \end{tikzpicture}
\caption{Polytope C}
\label{subfig:N=3 C}
\end{subfigure}}
\qquad
{\begin{subfigure}{0.35\textwidth}
\begin{tikzpicture}[scale=0.4] 
\WeylChamber
	\draw[very thick,fill=Crimson] (0,1.5) node[anchor=east]{$a$} -- 
	(1.299, 3.751) circle [fill,radius=0.08] node[black,anchor=west]{$c_2$} --
	(2.598, 6) circle [fill,radius=0.08] node[black,anchor=west]{$c_3$} -- 
	(3.897, 8.25) node[black,anchor=west]{$c_1$} --
	(1.299, 8.25) node[black,anchor=south]{$b$} --
	(0, 6) -- (0,1.5);
	\draw[thick] (2.598, 6) -- (1.299, 8.25);
	\draw[thick] (1.299, 3.751) -- (0,6); 
 \end{tikzpicture}
\caption{Polytope D} 
\label{subfig:N=3 D}
\end{subfigure}}


{\begin{subfigure}{0.35\textwidth}
\setcounter{subfigure}{4}
\begin{tikzpicture}[scale=0.4] 
\WeylChamber
	\draw[very thick,fill=Crimson] (0,4.5) node[anchor=east]{$a$} -- 
	(2.598, 9) node[black,anchor=south]{$c_3$} --
	(3.897, 6.7505) node[black,anchor=west]{$b$} --
	(1.299, 2.25) node[black,anchor=north]{$c_2$} -- (0,4.5);
	\draw[thick] (1.299, 6.75) -- (3.897, 6.7505);
	\draw [fill] (1.299, 6.75) circle [fill,radius=0.1] node[black,anchor=east]{$c_1$};
 \end{tikzpicture}
\caption{Polytope E}
\label{subfig:N=3 E}
\end{subfigure}}
\qquad
{\begin{subfigure}{0.35\textwidth}
\begin{tikzpicture}[scale=0.4] 
\WeylChamber
	\draw[very thick,fill=Crimson] (0,3) node[anchor=east]{$a$} -- 
	(3.118, 8.4) node[black,anchor=south]{$c_3$} --
	(4.936, 5.25) node[black,anchor=west]{$b$} -- 
	(2.857, 1.65) -- (1.299, .75) -- 
	(.7794, 1.65)  -- (0,3);
	\draw[thick] (1.299, 5.25) -- (4.936, 5.25);
	\draw[thick] (.7794, 1.65) -- (2.857, 1.65);
	\draw [fill] (.7794, 1.65) circle [radius=0.1] node[black,anchor=east]{$c_2$};
	\draw [fill] (1.299, 5.250) circle [radius=0.1] node[black,anchor=east]{$c_1$};
\end{tikzpicture}
\caption{Polytope F}
\label{subfig:N=3 F}
\end{subfigure}}

\vskip 0.5cm 

\begin{subfigure}{0.35\textwidth}
\begin{tikzpicture}[scale=0.4] 
\WeylChamber
	\draw[very thick, fill=Crimson] (0,1.5) node[anchor=east]{$a$} -- 
	(.7794, 2.85) circle [fill,radius=0.08] node[black,anchor=east]{$c_2$} --
	(1.819, 4.65) circle [fill,radius=0.08] node[black,anchor=east]{$c_1$} -- 
	(3.897, 8.25) node[black,anchor=south]{$c_3$} --
	(5.975, 4.65) node[black,anchor=west]{$b$} -- 
	(4.936, 2.850) -- (0.6494, 0.375) -- 
        (0,1.5);
	\draw[thick] (1.819, 4.65) -- (5.975, 4.65);
	\draw[thick] (.7794, 2.85) -- (4.936, 2.8);
 \end{tikzpicture}
\caption{Polytope G}
\label{subfig:N=3 G}
\end{subfigure}
\qquad
{\begin{subfigure}{0.35\textwidth}
\begin{tikzpicture}[scale=0.4] 
\WeylChamber
	\draw[very thick,fill=Crimson] (0,5.1) node[anchor=east]{$a$} -- 
	(2.078, 8.73) node[black,anchor=south]{$c_3$} --
	(4.676, 4.2) node[black,anchor=west]{$b$} --
	(3.377, 1.95) -- (2.208, 1.276) --
	(1.819, 1.95) -- (.5196, 4.22) -- (0,5.1);
	\draw[thick] (.5196, 4.22) -- (4.676, 4.2);
	\draw[thick] (3.377, 1.95) -- (1.819, 1.95);
	\draw (.5196, 4.22) circle [fill,radius=0.08] node[black,anchor=north]{$c_1\;$};
	\draw (1.819, 1.95) circle [fill,radius=0.08] node[black,anchor=east]{$c_2$};

 \end{tikzpicture}
\caption{Polytope H}
\label{subfig:N=3 H}
\end{subfigure}}
\caption{The generic momentum polytopes: refer to Fig.\,\ref{fig:parameter plane} for the notation.}
\label{fig:polytopes A-H}
\end{figure}
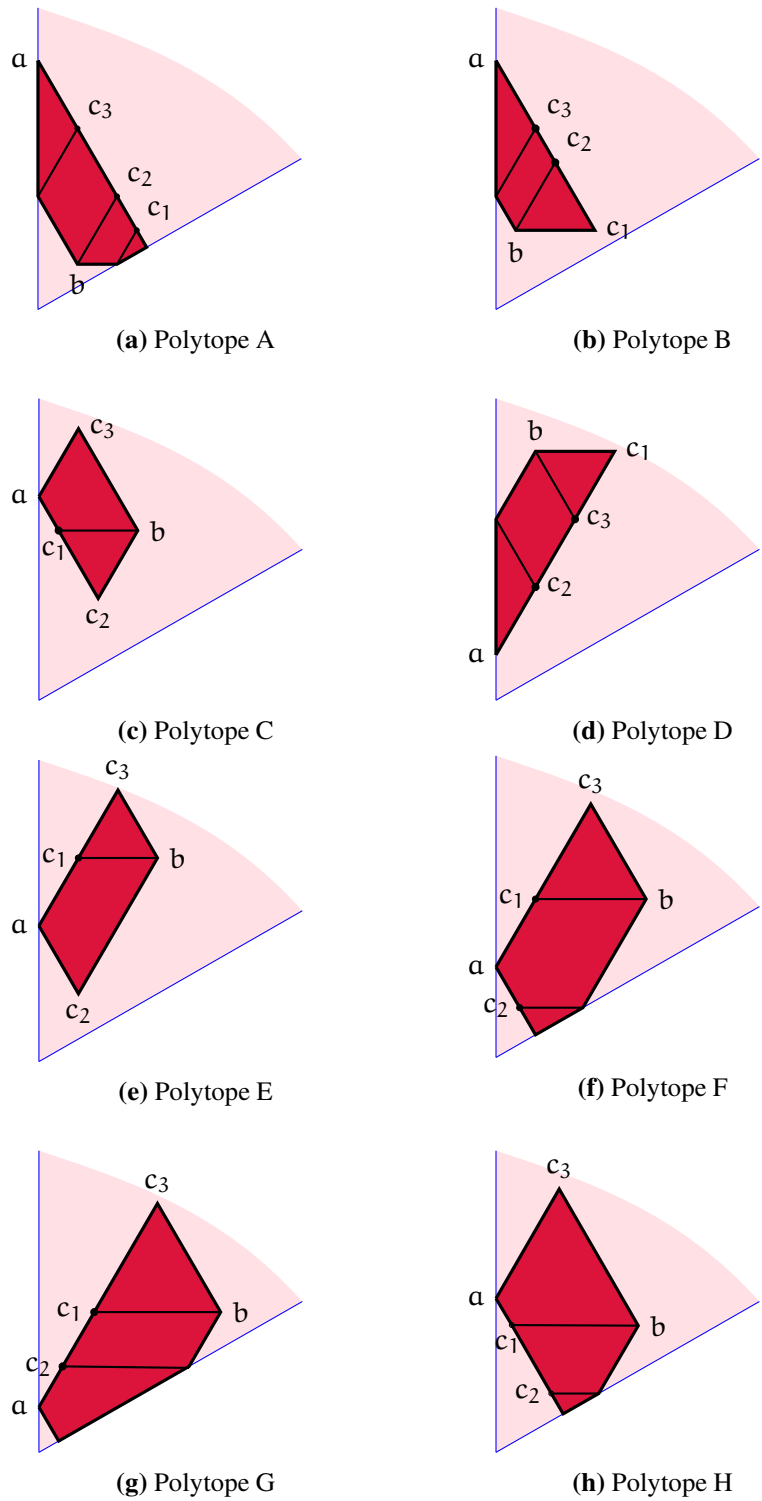

\begin{theorem} \label{thm:3-polytopes} 
The momentum polytopes of the SU(3) action on $\CP ^2\times\CP ^2\times\CP ^2$ with weighted symplectic form $\Omega=\Gamma_1\omega_{FS}\oplus\Gamma_2\omega_{FS}\oplus\Gamma_3\omega_{FS}$ (with $\Gamma_j\neq0$) fall into eight different types for which $\Gamma_i-\Gamma_j-\Gamma_k\neq0$, $\Gamma_i\pm\Gamma_j\neq0$,  where $i,j,k=1,2,3$, and $\Gamma_1+\Gamma_2+\Gamma_3>0$. An example of each type is shown in Figure\,\ref{fig:polytopes A-H}, while the 8 different regions of $\Gamma$-space are illustrated in Figure\,\ref{fig:parameter plane}.
\end{theorem}

\begin{remark}\label{rmk:change sign}
There are another 8 types of generic momentum polytope with $\Gamma_1+\Gamma_2+\Gamma_3<0$.  Note from the expression for $\J$ above, if the signs of all three $\Gamma_j$ are changed, then the sign of $\J$ changes. This implies that the original polytope and the new one are related by the involution $*$ described in Remark\,\ref{rmk:involution}. 
The cases with $\Gamma_1+\Gamma_2+\Gamma_3=0$ are transitional and described with other transitional cases further below---see Figure\,\ref{fig:sum=0}.
\end{remark}

The remainder of this section consists of a proof of this theorem.  See Remark\,\ref{rmk:types of polytope} for a discussion of the word `type'. Non-generic, or transitional, polytopes are discussed in Section\,\ref{sec:transitions}, and illustrated in diagrams at the end of the paper.

Recall from Sjamaar's theorem \ref{thm:Sjamaar} that if a point $\mu$ in the interior of the positive Weyl chamber is a vertex of the momentum polytope then there is an $m\in M$ with $\J(m)=\mu$ and stabilizer equal to a maximal torus.   We begin therefore with an analysis of the points with stabilizer equal to a maximal torus, which we choose to be the subgroup $\TT$ of diagonal matrices.  

The fixed points of $\TT$ are the 27 points $(e_i,e_j,e_k)\in M$, where each $e_\ell\in\{e_1,e_2,e_3\}$ and $e_1=[1:0:0],\,e_2=[0:1:0]$ and $e_3=[0:0:1]$.  The images in $\tt^*_+$ of these points are determined by the spectra $\sigma$ of the corresponding Hermitian matrix.  Extending the notation for $\CP^2\times\CP^2$, let
$$
a=\sigma\Bigl(J(e_1,e_1,e_1)\Bigr),\qquad
b=\sigma\Bigl(J(e_1,e_2,e_3)\Bigr),$$
$$c_1=\sigma\Bigl(J(e_2,e_1,e_1)\Bigr),\quad
c_2=\sigma\Bigl(J(e_1,e_2,e_1)\Bigr),\quad
c_3=\sigma\Bigl(J(e_1,e_1,e_2)\Bigr) .
$$
Here some caution is required: for example $(e_1,e_1,e_1)$ and $(e_2,e_2,e_2)$ and $(e_3,e_3,e_3)$ lie on the same $SU(3)$-orbit in $M$, and hence their values under $J$ lie on the same Weyl group orbit in $\tt^*$.  Only one of these will lie in the positive Weyl chamber (we denote this value $a$). On the other hand, their (unordered) spectra coincide, and so we consider spectra as sets. The same applies to say $(e_1,e_2,e_1)$---permuting the indices to, for example, $(e_2,e_3,e_2)$ will give points in the same orbit (but not the same orbit as $(e_2,e_1,e_1)$), so their unordered spectra are equal (but in general different to that of $(e_2,e_1,e_1)$).  

One finds 
\begin{equation}\label{eq:vertices}
\begin{array}{lcl}
a &=& \Bigl\{\tfrac23\left(\Gamma_1+\Gamma_2+\Gamma_3\right) , 
	-\tfrac13\left(\Gamma_1+\Gamma_2+\Gamma_3\right) ,
	-\tfrac13\left(\Gamma_1+\Gamma_2+\Gamma_3\right) \Bigr\},\\[6pt]
b &=&  \Bigl\{\tfrac13\left(2\Gamma_1-\Gamma_2-\Gamma_3\right) , 
	\tfrac13\left(-\Gamma_1+2\Gamma_2-\Gamma_3\right) ,
	\tfrac13\left(-\Gamma_1-\Gamma_2+2\Gamma_3\right) \Bigr\},\\[6pt]
c_1 &=&   \Bigl\{\tfrac13\left(2\Gamma_1-\Gamma_2-\Gamma_3\right) , 
	\tfrac13\left(-\Gamma_1+2\Gamma_2+2\Gamma_3\right) ,
	-\tfrac13\left(\Gamma_1+\Gamma_2+\Gamma_3\right) \Bigr\},\\[6pt]
c_2 &=&  \Bigl\{\tfrac13\left(-\Gamma_1+2\Gamma_2-\Gamma_3\right) , 
	\tfrac13\left(2\Gamma_1-\Gamma_2+2\Gamma_3\right) ,
	-\tfrac13\left(\Gamma_1+\Gamma_2+\Gamma_3\right) \Bigr\},\\[6pt]
c_3 &=& \Bigl\{\tfrac13\left(-\Gamma_1-\Gamma_2+2\Gamma_3\right) , 
	\tfrac13\left(2\Gamma_1+2\Gamma_2-\Gamma_3\right) ,
	-\tfrac13\left(\Gamma_1+\Gamma_2+\Gamma_3\right) \Bigr\}.
\end{array}
\end{equation}
In order to depict them in the positive Weyl chamber, each of these sets should be ordered by $\lambda_1\geq\lambda_2\geq\lambda_3$.  Note that point $a$ always lies on a wall of the Weyl chamber as two of the eigenvalues are equal.  The other points do not lie on a wall in general (when the 3 eigenvalues are distinct), however for some values of the weights the points can lie on the wall and these determine the transition cases. 

For example, $b$ lies on a wall if two of the weights coincide, while $c_1$ lies on a wall if $\Gamma_1=\Gamma_2+\Gamma_3$ or $\Gamma_2+\Gamma_3=0$ (or the degenerate the case $\Gamma_1=0$ which we exclude from disussions).  Similar possibilities occur for $c_2$ and $c_3$ with the indices of the $\Gamma_j$ permuted accordingly.  The set of possible degeneracies, up to permutations of the  indices, are listed in Table\,\ref{table:transitions}.

\begin{table}
$$
\begin{array}{r|l}
\text{condition} & \text{degeneracy} \\
\hline
\rule{0pt}{12pt}\Gamma_1=0 & a=c_1,\; b=c_2=c_3\\[4pt]
\Gamma_1=\Gamma_2 & b\in \text{Wall},\;c_2=c_3\\[4pt]
\Gamma_1+\Gamma_2=0 &  a=c_3\,(\in \text{Wall})\\[4pt]
\Gamma_1=\Gamma_2+\Gamma_3 &  c_1\in \text{Wall}\\[4pt]
\Gamma_1+\Gamma_2+\Gamma_3=0 &  a=0
\end{array}
$$
\centering
\caption{Transition values of $\Gamma_j$ ; similar transitions occur permuting the indices. `$x\in\text{Wall}$' means that the point $x$ belongs to a wall of the Weyl chamber. See Figure\,\ref{fig:parameter plane}; further details are shown in Section\,\ref{sec:transitions} and Figures\,\ref{fig:transition parameters}--\ref{fig:transitions1}.}
\label{table:transitions}
\end{table}

For some values of the symplectic weights $\Gamma_j$, the convex hull of these 5 points is equal to the momentum polytope.  But for others we need to determine the `local momentum cones', which are determined by the local images of the orbit momentum map near these points.

There are two procedures that can be used for drawing the different momentum polytopes.  One is starting with one we know (eg for 2 copies of $\CP^2$ by putting one of the $\gamma_j$ to 0) and then varying the weights and following the possible polytope, and the other is looking at the local momentum cones for each vertex. In this paper we use mostly the local momentum cones, with some continuity arguments, while in the thesis \cite{AS-thesis} the former approach is used more.

\subsection{Generic polytopes}
We now proceed to calculate the local momentum cones at each of the 5 vertices $a,b,c_1,c_2$ and $c_3$. To do this we need to calculate $\J_{N_1}$ for each. 
At each of the $\TT$-fixed points $m=(e_i,e_j,e_k)$ the tangent space at $m$ is given by
$$T_mM = T_{e_i}\CP^2\;\times\;T_{e_j}\CP^2\;\times\;T_{e_k}\CP^2,$$
and this (symplectic) decomposition is invariant under the  action of the maximal torus $\TT$; see Lemma\,\ref{lemma:weights} for the weights of this action. 

\medskip

\noindent\underline{\em Vertex $b$} \quad
Consider the weights at $b=\J(m)$ for $m=(e_1,e_2,e_3)$.  If we put, 
$$\xx=\begin{pmatrix}\begin{pmatrix}0\\v_1\\w_1\end{pmatrix},\begin{pmatrix}u_2\\0\\ w_2\end{pmatrix},\begin{pmatrix}u_3\\v_3\\0\end{pmatrix}\end{pmatrix} \in T_mM,$$
then, 
$$D\J_m (\xx) = \Gamma_1\begin{pmatrix}
0&\bar v_1&\bar w_1\\ v_1&0&0 \\ w_1&0&0
\end{pmatrix} +
\Gamma_2\begin{pmatrix}
0&u_2&0\\ \bar u_2&0& \bar w_2 \\ 0&w_2&0
\end{pmatrix} 
+\Gamma_3\begin{pmatrix}
0&0&u_3\\ 0&0&v_3 \\ \bar u_3&\bar v_3&0
\end{pmatrix} 
$$
Thus $\ker D\J_m$ consists of those $\xx$ satisfying
$$\Gamma_1\bar v_1+\Gamma_2u_2=0,\quad \Gamma_2\bar w_2+\Gamma_3v_3=0,\quad \Gamma_1w_1+\Gamma_3\bar u_3=0.$$
This defines a subspace of dimension 6. For this section, we assume $b$ is not contained in a wall of the Weyl chamber, and then this is in fact the symplectic slice: whenever $G_m=G_\mu$ one has $T_0=0=N_0$, and thus $N_1=\ker D\J_m$.  

To find $\J_{N_1}$ is simple: since $N_1$ is the sum of 3 distinct representations, with weights $\alpha_1,\alpha_2$ and $\alpha_3$ respectively, the momentum map is a sum of three terms (see \eqref{eq:Gamma_j mom map}). Thus, using $w_1,u_2,v_3$ to parametrize $N_1$ (with $v_1=-(\Gamma_2/\Gamma_1)\bar u_2$ etc.)
\begin{equation}\label{eq:JN1 for b}
  \J_{N_1}(w_1,u_2,v_3) = 
	\frac{\Gamma_3}{\Gamma_2}(\Gamma_2-\Gamma_3)|v_3|^2\alpha_1 \ + \  
	\frac{\Gamma_1}{\Gamma_3}(\Gamma_3-\Gamma_1)|w_1|^2\alpha_2 \ + \  
	\frac{\Gamma_2}{\Gamma_1}(\Gamma_1-\Gamma_2)|u_2|^2\alpha_3.
\end{equation}
This determines the momentum cone at $b$ depending on the signs of the coefficients; it turns out that provided the three weights are distinct, this is always a 120$^\circ$ cone (if 2 of the weights coincide it becomes a 60$^\circ$ cone, but in that case $b$ is contained in a wall of the Weyl chamber---see further below for this case).  For example, if $\Gamma_1=4,\Gamma_2=2,\Gamma_3=-1$ (which lies in region C in Figure\,\ref{fig:parameter plane}) then 
$b=(7,1,-8) \in\tt^*_+$ and 
$$\J_{N_1}(u_2,v_3,w_1) =   -\tfrac32|v_3|^2\alpha_1 + 20\,|w_1|^2\alpha_2 + |u_2|^2\alpha_3,$$
whose image is precisely the cone at $b$ shown in Figure\,\ref{subfig:N=3 C} (see Figure\,\ref{fig:SU3rootdiagram} for the definition of the $\alpha_j$).

This expression is only the local momentum cone at $b$ provided $\J(m)\in\tt^*_+$~; if that is not the case then the calculation needs repeating for whichever of the 6 equivalent points does map to $b$. For example, if the weights $\Gamma_j$ are such that $b=\J(m)$ for $m=(e_2,e_1,e_3)$ then the expression for the local momentum cone (i.e., for $\J_{N_1}$) is like the one above, but with the roots permuted by the appropriate element of the Weyl group; thus, in that case,
$$\J_{N_1}(w_1,v_2,u_3) = 
	\frac{\Gamma_1}{\Gamma_3}(\Gamma_1-\Gamma_3)|w_1|^2\alpha_1 \ + \  
	\frac{\Gamma_3}{\Gamma_2}(\Gamma_3-\Gamma_2)|u_3|^2\alpha_2 \ + \  
	\frac{\Gamma_2}{\Gamma_1}(\Gamma_2-\Gamma_1)|v_2|^2\alpha_3.
$$
The choice we have made, that $\Gamma_1\geq\Gamma_2\geq\Gamma_3$, indeed ensures $\J(e_1,e_2,e_3)\in\tt^*_+$.

\medskip

\noindent\underline{\em Vertices $c_j$} \quad 
The calculations for $c_1,c_2$ and $c_3$ are very similar.  For example, with $m=(e_1,e_2,e_2)$ the elements of $T_mM$ can be written
$$\xx=\begin{pmatrix}\begin{pmatrix}0\\v_1\\w_1\end{pmatrix},\begin{pmatrix}u_2\\0\\ w_2\end{pmatrix},\begin{pmatrix}u_3\\0\\w_3\end{pmatrix}\end{pmatrix} \in T_mM.$$
Then
$$\xx\in\ker D\J_m \Longleftrightarrow w_1=(\Gamma_2w_2+\Gamma_3w_3)= (\Gamma_1\overline{v_1}+\Gamma_2u_2+\Gamma_3u_3)=0.$$ 
After similar calculations for other $c_j$, one obtains the following expressions for the slice momentum map.

\noindent$\bullet$ For $c_1$ using $m=(e_1,e_2,e_2)$,
\begin{equation}\label{eq:JN1 for c1}
J_{N_1} =
	 -\frac{\Gamma_2}{\Gamma_3}(\Gamma_2+\Gamma_3)|w_2|^2\alpha_1+ 
	 \left(\frac{\Gamma_1}{\Gamma_2}(\Gamma_1-\Gamma_2)|v_1|^2+
	 \frac{\Gamma_1\Gamma_3}{\Gamma_2}(u_3v_1+\overline{u_3v_1})+
	 \frac{\Gamma_3}{\Gamma_2}(\Gamma_3+\Gamma_2)|u_3|^2\right)\alpha_3.
\end{equation}
Lemma\,\ref{lemma:posdef} below shows that the coefficient of $\alpha_3$ is definite if and only if 
$$\Gamma_1\Gamma_2\Gamma_3(\Gamma_1-\Gamma_2-\Gamma_3)>0.$$ 
If this inequality is satisfied then $c_1$ lies at a vertex of the momentum polytope; if, on the other hand, it the expression is negative, then $c_1$ lies on an edge of the polytope (parallel to $\alpha_3$). The given inequality is satisfied only in regions B and D.  

\noindent$\bullet$ For $c_2$ using $m=(e_2,e_1,e_2)$, one obtains
\begin{equation}\label{eq:JN1 for c2}
J_{N_1} = 
-	\frac{\Gamma_3}{\Gamma_1}(\Gamma_1+\Gamma_3)|w_3|^2\alpha_1 + 
	\left(
	\frac{\Gamma_2}{\Gamma_3}(\Gamma_2-\Gamma_3)|v_2|^2
	+ 
	\frac{\Gamma_1\Gamma_2}{\Gamma_3}(u_1v_2+\overline{u_1v_2})
	+
	\frac{\Gamma_1}{\Gamma_3}(\Gamma_1+\Gamma_3)|u_1|^2
	\right)\alpha_3.
\end{equation}

\noindent$\bullet$ Finally, for $c_3$ using $m=(e_2,e_2,e_1)$,
\begin{equation}\label{eq:JN1 for c3}
J_{N_1} = 
	-\frac{\Gamma_1}{\Gamma_2}(\Gamma_1+\Gamma_2)|w_1|^2\alpha_1 + 
	\left(\frac{\Gamma_3}{\Gamma_1}(\Gamma_3-\Gamma_1)|v_3|^2 + 
	\frac{\Gamma_2\Gamma_3}{\Gamma_1}(u_2v_3+\overline{u_2v_3}) +
	\frac{\Gamma_2}{\Gamma_1}(\Gamma_1+\Gamma_2)|u_2|^2\right)\alpha_3.
\end{equation}
Similar conditions on the $\Gamma_j$ based on Lemma\,\ref{lemma:posdef} ensure the coefficients of $\alpha_3$ are definite or not.  

As with $b$, if $\J(m)$ fails to belong to the positive Weyl chamber, then the appropriate element of the Weyl group should be applied to the roots.  It turns out that $\J(e_1,e_2,e_2)\in\tt^*_+$ if and only if $\Gamma_1+\Gamma_2\geq\Gamma_3\geq0$.  Notice that while $\J_{N_1}$ for $b$ has all 3 roots appearing, for $c_1,c_2$ and $c_3$ it only has two distinct roots, which explains why only 2 lines (or half-lines) pass through those points in the figures.

\begin{figure}
	\centering
{\begin{subfigure}{0.4\textwidth}
\begin{tikzpicture}[scale=0.5] 
\WeylChamber
	\draw[thick,->] (0,1.5) node[anchor=east]{$a$} -- (0.6,2.539);
	\draw[thick, ->] 	(1.299, 3.751) node[black,anchor=west]{$c_2$} -- (0.7490,2.797);
	\draw[thick, ->] 	(1.299, 3.751)  -- (2.049,5.049);
	\draw[thick, ->] 	(1.299, 3.751)  -- (0.549,5.049);
	\draw (2.598, 6) circle [fill,radius=0.08] node[black,anchor=west]{$c_3$} ;
	\draw[thick, ->] 	(3.897, 8.25) node[black,anchor=west]{$c_1$} -- (2.7, 8.25);
	\draw[thick, ->] 	(3.897, 8.25) -- (3.147,6.951);
	\draw[thick, ->] 	(1.299, 8.25) node[black,anchor=south]{$b$} -- (2.5, 8.25);
	\draw[thick, ->] 	(1.299, 8.25) -- (0.5490,6.951);
	\draw[thick, ->] 	(1.299, 8.25) -- (2.049,6.951);
 \end{tikzpicture}
\caption{Weights for polytope D\\(weights at $c_3$ are translations of $c_2$)} 
\label{subfig:D-weights}
\end{subfigure}}
	\qquad\quad
\begin{subfigure}{0.4\textwidth}
\begin{tikzpicture}[scale=0.5] 
\WeylChamber
	\draw[thick, ->] (0,1.5) node[anchor=east]{$a$} -- (.50, 2.366);
	\draw[thick, ->] (0,1.5)  -- (0.5,0.634);
	\draw[thick, ->] 	(1.819, 4.65) node[black,anchor=east]{$c_1$} -- (1.069,3.351) ;
	\draw[thick, ->] 	(1.819, 4.65)  -- (2.569,5.949);
	\draw[thick, ->] 	(1.819, 4.65)  -- (3.319, 4.65);
	\draw[thick, ->] 	(.7794, 2.85) circle [fill,radius=0.08] node[black,anchor=east]{$c_2$} ;
	\draw[thick, ->] (3.897, 8.25) node[black,anchor=south]{$c_3$} -- (4.647,6.951);
	\draw[thick, ->] (3.897, 8.25)  -- (3.147,6.951);
	\draw[thick, ->] (5.975, 4.65) node[black,anchor=west]{$b$} -- (5.225,5.95);
	\draw[thick, ->] (5.975, 4.65)  -- (5.225,3.351); 
	\draw[thick, ->] (5.975, 4.65)  -- (4.475, 4.65); 
 	\draw (0,2) node {\hbox to 0pt{\hss%
		\begin{tikzpicture}[scale=0.18] 
			\draw [blue,->] (0,0) -- (4,0) node[anchor=north] {$\alpha_1$};
			\draw [rotate=-120,blue,->] (0,0) -- (4,0) node[anchor=north] {$\alpha_2$};
			\draw [rotate=120,blue,->] (0,0) -- (4,0) node[anchor=south] {$\alpha_3$};
		 \end{tikzpicture}\hskip8mm}};
 \end{tikzpicture}
\caption{Weights for polytope G\\(weights at $c_2$ are translations of $c_1$)} 
\label{subfig:G-weights}
\end{subfigure}
\caption{Examples showing weights at the fixed points}
\end{figure}
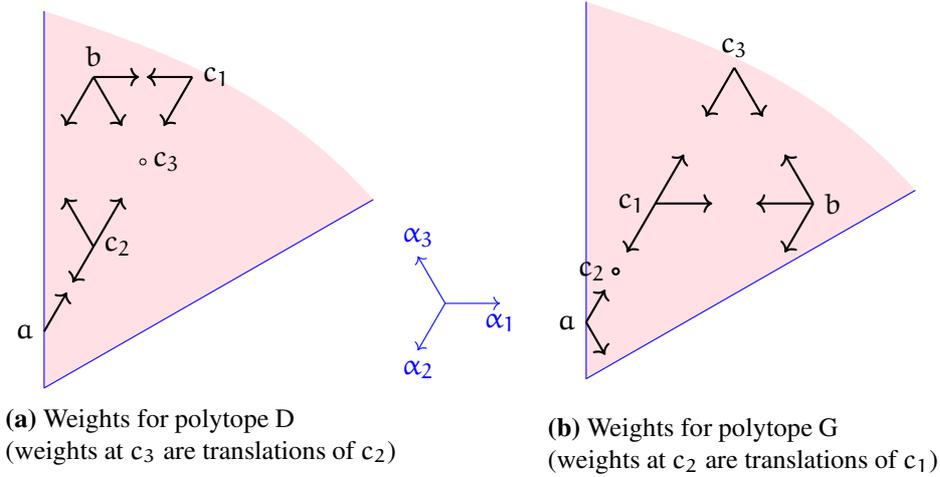

\begin{lemma}\label{lemma:posdef}
The real quadratic form $q:\C^2\to\R$ given by
$$q(u,v) = A|u|^2+B(u\overline{v}+\overline{u}v)+C|v|^2$$
is definite if and only if $AC>B^2$.
\end{lemma}

The proof is a simple calculation.

\medskip

\noindent\underline{\em Vertex $a$:} There remains to consider the point $a$.  Let $m$ be a point on the diagonal, and to be specific we take $m=(e_1,e_1,e_1)$.  The stabilizer of $m$ is now
$$G_m=\left\{\begin{pmatrix}
(\det A)^{-1}&0 \cr 0&  A
\end{pmatrix} \, \middle| \, A\in U(2)\right\} \simeq U(2).$$
We assume $\sum_j\Gamma_j\neq0$, in which case $a=\J(m)\neq0$ (see \eqref{eq:vertices}); the other case is covered by a continuity argument discussed below.   
We now have $\J_{N_1}:N_1\to\u(2)^*\subset\su(3)^*$. 

Since again $G_\mu=G_m$, it follows that $T_0=N_0=0$ and $N_1=\ker DJ_m$, which is of dimension 8.  Thus,
$$\xx\in N_1 \Longleftrightarrow \sum\Gamma_jv_j=\sum\Gamma_jw_j=0.$$
Before solving for $v_3$ and $w_3$ one finds
\begin{equation}\label{eq:JN1 for a}
\J_{N_1} = 
\begin{pmatrix}
-\sum_j\Gamma_j(|v_j|^2+|w_j|^2) & 0 & 0\\
0 & \sum_j\Gamma_j|v_j|^2 & \sum_j\Gamma_j \overline{v_j}w_j \\
0 & \sum_j\Gamma_jv_j\overline{w_j} &\sum_j\Gamma_j|w_j|^2 
\end{pmatrix}.
\end{equation}
With the usual Cartan subalgebra of diagonal matrices, there are two types of weight vector for the $\TT^2$ action on $T_mM$, vectors with $w_j=0$ (of weight $-\alpha_3$) and those with $v_j=0$ (of weight $\alpha_2$).  Now $N_1$ consists of 4-dimension's worth of each.  Eliminating $v_3$ and $w_3$ from the expression above, one finds, on the $\alpha_3$-weight space,
\begin{equation}\label{eq:JN1 for a-alpha3}
\J_{N_1}(v_1,v_2,0,0) = 
	\left(\frac{\Gamma_1}{\Gamma_3}(\Gamma_1+\Gamma_3)|v_1|^2 + 
	\frac{\Gamma_1\Gamma_2}{\Gamma_3}(v_1\overline{v_2}+\overline{v_1}v_2) +
	\frac{\Gamma_2}{\Gamma_3}(\Gamma_2+\Gamma_3)|v_2|^2\right)(-\alpha_3),
\end{equation}
and on the $\alpha_2$-weight space,
\begin{equation}\label{eq:JN1 for a-alpha2}
\J_{N_1}(0,0,w_1,w_2) = 
	\left(\frac{\Gamma_1}{\Gamma_3}(\Gamma_1+\Gamma_3)|w_1|^2 + 
	\frac{\Gamma_1\Gamma_2}{\Gamma_3}(w_1\overline{w_2}+\overline{w_1}w_2) +
	\frac{\Gamma_2}{\Gamma_3}(\Gamma_2+\Gamma_3)|w_2|^2\right)\alpha_2.
\end{equation}
For both of these, using the notation of Lemma\,\ref{lemma:posdef}, one finds the discriminant
$$AC-B^2= \frac{\Gamma_1\Gamma_2}{\Gamma_3}(\Gamma_1+\Gamma_2+\Gamma_3).\eqno({\bf D})$$
This expression \D\ is positive in regions A, B and D only.  Consider the different cases:

\smallskip 

\noindent\underline{\D${}<0$}\,: In this case the quadratic coefficients in \eqref{eq:JN1 for a-alpha3} and \eqref{eq:JN1 for a-alpha2} are indefinite, and the image of the momentum map at $a$ contains lines in the root directions $\pm\alpha_2,\pm\alpha_3$, and in the positive Weyl chamber this gives lines in the directions $-\alpha_2$ and $-\alpha_3$.  See for example Figure\,\ref{subfig:G-weights}.  The convexity theorem implies that the infinitesimal momentum cone at $a$ in this case contains the region between these two directions, but does not tell us if it is equal to it (see further below). 

\smallskip 

\noindent\underline{\D${}>0$}\,:  Here there are two possibilities: the quadratic coefficients in \eqref{eq:JN1 for a-alpha3} and \eqref{eq:JN1 for a-alpha2} are either \emph{positive} or \emph{negative} definite.  Suppose they are positive definite;  then the image of $J_{N_1}$ contains the directions $\alpha_2$ and $-\alpha_3$. However, from $a$, the direction $\alpha_2$ does not lie in the positive Weyl chamber, and applying the Weyl-group reflection fixing $a$ sends $\alpha_2$ to $-\alpha_3$. Thus all we know is that the image of $\J$ in a neighbourhood of $a$ contains a line in the direction of $-\alpha_3$ (see Figures\,\ref{subfig:N=3 A} and \ref{subfig:N=3 B}).  Similarly, if they are negative definite, the image contains a line in the direction of $-\alpha_2$ (as in Figure\,\ref{subfig:N=3 D}).

\medskip

\def\cV{\mathcal{V}}
\noindent\underline{\emph{Conclusion \& construction of generic polytopes:}}  
For $\Gamma=(\Gamma_1,\Gamma_2,\Gamma_3)$ in each of the regions A,B,\dots,H of the diagram in Figure\,\ref{fig:parameter plane}, one plots the five points $a$ (on the wall) and $b,c_1,c_2,c_3$ in the interior of the positive Weyl chamber. From each point, one can plot the local momentum cone. The  theorem of Sjamaar (see Theorem\,\ref{thm:Sjamaar}(2) above) states,
$$\Delta(M) = \bigcap_{m\in\Phi^{-1}(\tt^*_+)} \Delta_m$$
where $\Delta_m$ is the local momentum cone at $\Phi(m)$ (independent of $m$ in the fibre), as defined above.  As can be seen from the figures, the points $a,b$ and the $c_j$ only account for some of the vertices (there may be others on the boundary of the Weyl chamber). However, if we put $\cV=\{a,b,c_1,c_2,c_3\}$, it follows from Sjamaar's theorem that
\begin{equation}\label{eq:vertex LMCs}
\Delta(M) \subset \bigcap_{\mu\in \cV} \Delta_\mu.
\end{equation}
In each region except G, the information from $b,c_1,c_2,c_3$ suffices. For example, for $\Gamma$ in region D, refer to the weights shown in Figure\,\ref{subfig:D-weights}.  Starting from the point $c_1$, the weights dictate a line from $c_1$ to $b$, and from $c_1$ to $c_3$ to $c_2$ and thence to $a$. From $b$ the weight in the direction $\alpha_2$ leads to the wall.  The convex hull of this set is the unique set satisfying the inclusion \eqref{eq:vertex LMCs} and in addition containing the infinitesimal momentum cones (see Definition\,\ref{def:LMC} and Theorem\,\ref{thm:Sjamaar}).  

\medskip

\begin{figure}[t]
\centering
{\begin{subfigure}{0.2\textwidth}
\begin{tikzpicture}[scale=1]
 \draw[blue,fill=Pink!50] (0,2) -- (0,0) -- (1.732, 1);
	\draw[thick, fill=OrangeRed] (.4330, 1.75) -- (0,1) node[anchor=east]{$a$} -- (.4330, .25) -- (1.732, 1);
	\draw (0,0) node[anchor=east]{$0$};
	\fill[black] (0.2165,0.625) circle (0.05)  node[anchor=north]{$g$};
\end{tikzpicture}
\caption{\hskip1cm\ }
\label{subfig:point G(a)}
\end{subfigure}}
	\qquad
\begin{subfigure}{0.2\textwidth}
\begin{tikzpicture}[scale=1]
 \draw[blue,fill=Pink!50] (0,2) -- (0,0) -- (1.732, 1);
	\draw[thick, fill=OrangeRed] (.4330, 1.75) -- (0,1) node[anchor=east]{$a$} -- (0,0.5) --(.2165, .125) -- (1.732, 1);
	\draw[thick] (0,1) -- (.4330, .25);
	\draw (0,0) node[anchor=east]{$0$};
	\fill[black] (0.2165,0.625) circle (0.05)  node[anchor=south west]{$g$};
\end{tikzpicture}
\caption{\hskip1cm\ }
\label{subfig:point G(b)}
\end{subfigure}
	\qquad
\begin{subfigure}{0.2\textwidth}
\begin{tikzpicture}[scale=1]
 \draw[blue,fill=Pink!50] (0,2) -- (0,0) -- (1.732, 1);
	\draw[thick, fill=OrangeRed] (.4330, 1.75) -- (0,1) node[anchor=east]{$a$} -- (0,0) -- (1.732, 1);
	\draw (0,0) node[anchor=east]{$0$};
	\draw[thick] (0,1) -- (.4330, .25);
	\fill[black] (0.2165,0.625) circle (0.05) node[anchor=south west]{$g$};;
\end{tikzpicture}
\caption{\hskip1cm\ }
\label{subfig:point G(c)}
\end{subfigure}
\caption{Three possibilites for the lower part of polytope G compatible with local information at vertices $b,c_1,c_2,c_3$---version (a) is the correct one as shown by considering the local momentum cone at $g$.}
\label{fig:poly G possibilities}
\end{figure}
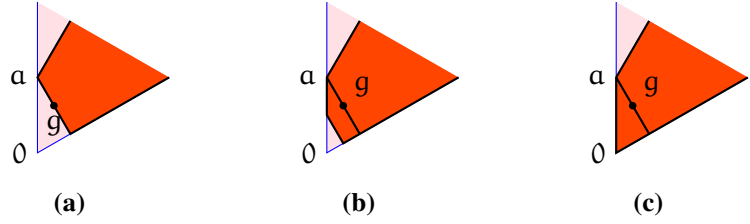

\noindent\underline{Region G:} The argument above suffices for all the generic polytopes except those of region G; the three diagrams in Figure\,\ref{fig:poly G possibilities} are all compatible with the data at vertices $b,c_1,c_2,c_3$, and we need to consider in greater detail the local momentum cone at $a$.  This region G is defined by the inequalities 
$-\Gamma_3>\Gamma_1>\Gamma_2>\Gamma_3$ and $\Gamma_1+\Gamma_2+\Gamma_3>0$ (see Figure\,\ref{fig:parameter plane}), and hence expression \D\ is negative. We need to determine in particular whether, at a point $g$ of the line in the direction $-\alpha_3$ (see Figure\,\ref{fig:poly G possibilities}), the infinitesimal momentum cone is the germ of a half space (above and to the right of the line, as in Figure\,\ref{subfig:point G(a)}) or of the full space (as in Figures\,\ref{subfig:point G(b)}, \ref{subfig:point G(c)}).

To accomplish this, consider the point $m'$ with $v_1=1,v_2=w_1=w_2=0$ in the symplectic slice $N_1$ at $m=(e_1,e_1,e_1)$; we have $\J_{N_1}(m') = A(-\alpha_3)$, and put $g=a+\J_{N_1}(m')$.  Here $A=(\Gamma_1+\Gamma_3)\Gamma_1/\Gamma_3<0$ in region $G$. (One could replace $v_1=1$ with $v_1=\epsilon$ to ensure $g$ is in the positive Weyl chamber, but the calculation is identical save for a factor of $\epsilon^2$.)  Let us calculate the infinitesimal momentum cone $\delta'$ at this point. For the $U(2)$ action on $N_1$, the point $m'$ has stabilizer $U(1)$ generated by $\diag[i,i,-2i]\in\alpha_3^\circ$ (cf.\ Lemma\,\ref{eq:bifurcation lemma}). Moreover $G_g=\TT^2$.  Considering the $U(2)$-action on $N_1$ gives rise to a  Witt-Artin decomposition at $m'$ given by 
$$T_{m'}(N_1) = T_0'\oplus T_1'\oplus N_1'\oplus N_0'$$
with $\dim T_0'=\dim N_0'=1,\dim T_1'=2$ leaving $\dim N_1'=4$. Now, a local calculation shows that
$$N_1'=\{(v_1,v_2,w_1,w_2)\in N_1 \mid Av_1+Bv_2=Aw_1+Bw_2=0\},$$
where $A,B$ are as before the coefficients in \eqref{eq:JN1 for a-alpha3} above.  This space can be parametrized by $v_1=v,w_1=w$ and hence $v_2=-(A/B)v,\,w_2=-(A/B)w$. The image of $N_0'$ under the momentum map is the line along the root direction $\pm\alpha_3$. Moreover so is the image of $(v,0)$. Indeed, the $U(2)$-momentum map on $N_1'$ is
$$\J_{N_1'}(v,w) =\frac{A}{B^2}(AC-B^2) \begin{pmatrix}
-|v|^2-|w|^2&0&0\cr 0 & |v|^2&\overline{v}w\cr 0&v\overline{w}&|w|^2
\end{pmatrix}.
$$
This lies in $\tt^*$ if and only if $\overline{v}w=0$. The case $(v,0)$ has been mentioned, leaving the case $(0,w)$:
$$\J_{N_1'}(0,w) = \frac{A}{B^2}(AC-B^2)|w|^2\diag[-1,0,1] = \frac{A}{B^2}(AC-B^2)|w|^2\alpha_2.$$
Now, in region G, the coefficient $(A/B^2)(AC-B^2)<0$, and since the infinitesimal momentum cone $\delta'$ is independent of the point in the fibre (Theorem\,\ref{thm:Sjamaar}), this shows that the infinitesimal momentum cone is indeed the half space as claimed. 

With this the proof of Theorem\,\ref{thm:3-polytopes} is concluded. \hfill$\Box$


\begin{figure}[t]  
	\centering


  \begin{tikzpicture}[scale=1.5]   

	\draw [thick, dashed] (1.732,-1) -- (0.866,0.5);
	\draw [thick,dashed] (1.732,-1) -- (4.5,-1);
	
	\draw [blue] (0.433, -0.25) -- (0.866,0.5) ;
	\draw [blue] (1.732,-1) -- (3.464,2);

	\draw [blue] (0.866,0.5) -- (4.5,0.5);
	\draw [blue] (3.464,2) -- (4.5,2);
	
	\draw[line width=2.5pt,red] (0,0) -- (4.5,2.598);
	\draw[line width=2.5pt,red] (0,0) -- (4.5,-2.598);
	
	\draw[<-] (0.217,-0.125) -- (-0.5,-0.4) node [anchor=east]{AA};

\def\red{\color{red}}
	\draw [fill] (0,0) circle [radius=0.03] node [anchor=east]{AAA};
	\draw [fill] (0.58,0) circle [radius=0.03] node [anchor=west]{AB};
	\draw [fill] (0.433, -0.25) circle [radius=0.03] node [anchor=north east]{AAB};
	\draw [fill] (0.433,0.25) circle [radius=0.03] node [anchor=south east]{AA};
	\draw [fill] (0.217,-0.125) circle [radius=0.03]; 
	\draw [fill] (1.082,-.625) circle [radius=0.03] node [anchor=north east]{BB};
	\draw [fill] (2.598,0.5) circle [radius=0.03] node [anchor=north west]{\!\!CEFH};
	\draw [fill] (2.165,-0.25) circle [radius=0.03] node [anchor=west]{CE};
	\draw [fill] (3.464,2) circle [radius=0.03] node [anchor=south east]{FGH};
	\draw [fill] (3.031 ,1.25 ) circle [radius=0.03] node [anchor=west]{FH};
	\draw [fill] (3.982,2.299) circle [radius=0.03] node [anchor=south east]{GG};
	\draw [fill] (2.165,1.25) circle [radius=0.03] node [anchor=south east]{HH};
	\draw [fill] (3.464,-2) circle [radius=0.03] node [anchor=north east]{DD};
	\draw [fill] (1.8,0.5) circle [radius=0.03] node [anchor=north]{CH};
	\draw [fill] (3.982,0.5) circle [radius=0.03] node [anchor=north]{EF};
	\draw [fill] (3.982,2) circle [radius=0.03] node [anchor=north]{FG};

	\draw [blue,->] (4.5,-2.598) -- (4.95,-2.858) node [anchor=south west]{DD$_0$};
	\draw [blue,->] (4.5,-1.6) -- (4.95,-1.7) node [anchor=west]{D$_0$};
	\draw [blue,->] (4.5,2.2) -- (4.95,2.3) node [anchor=west]{G$_0$};
	\draw [blue,->] (4.5,2.598) -- (4.95,2.858) node [anchor=west]{GG$_0$};
  \end{tikzpicture}

\vskip 5mm
\caption{This shows the labels of all 20 transition polytopes with $\Gamma_j\neq0$. Compare with Fig.\,\ref{fig:parameter plane}. The transitions denoted D$_0$, DD$_0$, G$_0$ and GG$_0$ arise `at infinity' in this diagram, and refer to points with $\Gamma_1+\Gamma_2+\Gamma_3=0$; the polytopes are illustrated in Figure\,\ref{fig:sum=0}. The transition between D$_0$ and G$_0$ occurs when $\Gamma_2=\Gamma_1+\Gamma_3=0$.}
\label{fig:transition parameters}
\end{figure}
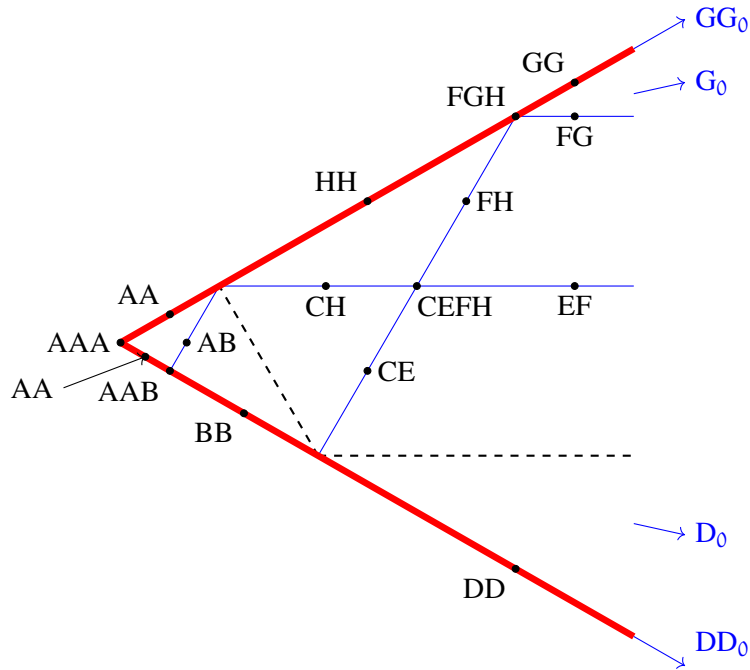

\subsection{Transition polytopes}  \label{sec:transitions}
\begin{figure}
	\centering

{\begin{subfigure}{0.3\textwidth}
\begin{tikzpicture}[scale=0.4] 
\WeylChamber
	\draw[very thick,fill=Crimson] (0, 9) -- (3.248, 3.375) -- (3.898, 2.250) --
	 (1.948, 1.125) --  (0, 4.5) -- (0, 9);
	\draw[thick] (1.948, 1.125) -- (3.248, 3.375);
	\draw[thick] (1.299, 6.750) -- (0,4.5);
	 \draw [fill] (0, 9) node[anchor=east]{$a$};
	 \draw [fill] (1.948, 1.125) node[anchor=north]{$b$};
	 \draw [fill] (3.248, 3.375) circle [fill,radius=0.1] node[anchor=south west]{$c_1=c_2$};
	 \draw [fill] (1.299, 6.75) circle [fill,radius=0.1] node[anchor=south west]{$c_3$};
 \end{tikzpicture}
\caption{Polytope AA  \\ with $\Gamma_1=\Gamma_2$}
\label{subfig:N=3 AA}
\end{subfigure}}
	\quad
{\begin{subfigure}{0.3\textwidth}
\begin{tikzpicture}[scale=0.4] 
\WeylChamber
	\draw[very thick,fill=Crimson] (0, 9) -- (0,0) -- (3.897, 2.25) --
	 (2.598, 4.5) -- (0, 9);
	\draw[thick] (0,0) -- (2.598, 4.5);
	 \draw [fill] (0, 9) node[anchor=east]{$a$};
	 \draw [fill] (0,0) node[anchor=east]{$b$};
	 \draw [fill] (2.598, 4.5) circle [fill,radius=0.1] node[anchor=west]{$c_1=c_2=c_3$};
 \end{tikzpicture}
\caption{Polytope AAA with \\ $\Gamma_1=\Gamma_2=\Gamma_3>0$}
\label{subfig:N=3 AAA}
\end{subfigure}}
	\quad
{\begin{subfigure}{0.3\textwidth}
\begin{tikzpicture}[scale=0.4] 
\WeylChamber
	\draw[very thick,fill=Crimson] (0, 9) -- (3.637, 2.7) -- (3.896, 2.25) --
	 (3.118, 1.8) --(0, 1.8) -- (0, 9);
	\draw[thick] (0,1.8) -- (2.078, 5.4);
	\draw[thick] (3.637, 2.70) -- (3.118, 1.8);
	 \draw [fill] (0, 9) node[anchor=east]{$a$};
	 \draw [fill] (0,1.8) node[anchor=east]{$b$};
	 \draw [fill] (3.637, 2.7) circle [fill,radius=0.1] node[anchor=south west]{$c_1$};
	 \draw [fill] (2.078, 5.4) circle [fill,radius=0.1] node[anchor=west]{$c_2=c_3$};
 \end{tikzpicture}
\caption{Polytope AA  \\ with $\Gamma_2=\Gamma_3$}
\label{subfig:N=3 AA'}
\end{subfigure}}

\vskip 5mm
\caption{The transition polytopes with repeated weights around region A.}
\label{fig:transitions2}
\end{figure}
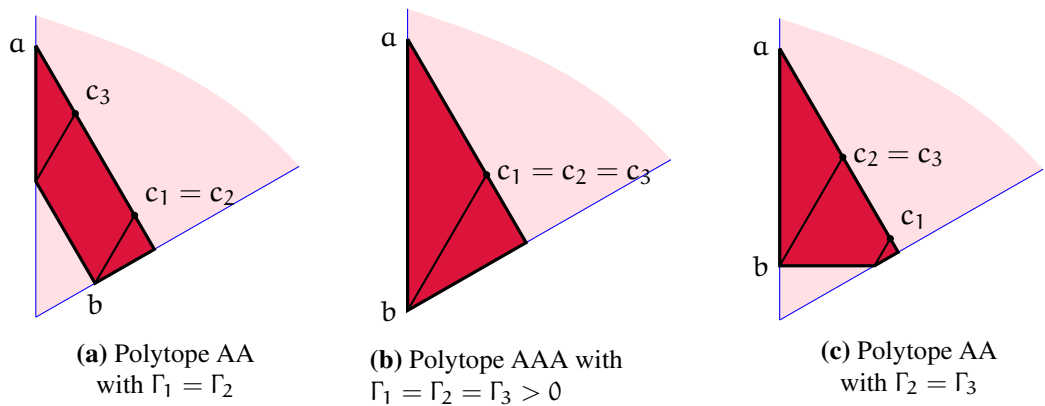

\begin{figure}[p] 
\centering

\centering

{\begin{subfigure}{0.3\textwidth}
\begin{tikzpicture}[scale=0.38] 
\WeylChamber
	\draw[very thick,fill=Crimson] (0, 9) -- (0,2.25) -- (3.897, 2.25) -- (0, 9);
	\draw[thick] (0., 2.25) -- (1.948, 5.625);
	 \draw [fill] (0, 9) node[anchor=east]{$a$};
	 \draw [fill] (0, 2.25) node[anchor=east]{$b$};
	 \draw  (3.897, 2.25) node[anchor=north]{$c_1$};
	 \draw [fill] (1.948, 5.625) circle [fill,radius=0.1] node[anchor=west]{$c_2=c_3$};
 \end{tikzpicture}
\caption{Polytope AAB}
\label{subfig:N=3 AAB}
\end{subfigure}}
	\quad
{\begin{subfigure}{0.3\textwidth}
\begin{tikzpicture}[scale=0.38] 
\WeylChamber
	\draw[very thick,fill=Crimson] (0, 9) -- (3.637, 2.7) -- (0, 2.7) -- (0, 9);
	\draw[thick] (0., 2.7) -- (1.819, 5.85);
	 \draw [fill] (0, 9) node[anchor=east]{$a$};
	 \draw [fill] (0, 2.7) node[anchor=east]{$b$};
	 \draw (3.637, 2.7)  node[anchor=north]{$c_1$};
	 \draw [fill] (1.819, 5.85) circle [fill,radius=0.1] node[anchor=west]{$c_2=c_3$};
 \end{tikzpicture}
\caption{Polytope BB}
\label{subfig:N=3 BB}
\end{subfigure}}
	\quad
{\begin{subfigure}{0.3\textwidth}
\begin{tikzpicture}[scale=0.38] 
\WeylChamber
	\draw[very thick,fill=Crimson] (0, 0.6) -- (0, 7.8) -- 
	 (4.157, 7.8)  -- (0, 0.6);
	\draw[thick] (0, 7.8) -- (2.078, 4.2);
	 \draw [fill] (0, 0.6) node[anchor=east]{\hskip5mm$a$};
	 \draw [fill] (0, 7.8) node[anchor=east]{$b$};
	 \draw (4.157, 7.8) node[anchor=west]{$c_1$};
	 \draw [fill] (2.078, 4.2) circle [fill,radius=0.1] node[anchor=west]{$c_2=c_3$};
 \end{tikzpicture}
\caption{Polytope DD}
\label{subfig:N=3 DD}
\end{subfigure}}

\vskip 5mm

{\begin{subfigure}{0.3\textwidth}
\begin{tikzpicture}[scale=0.38] 
\WeylChamber
	\draw[very thick,fill=Crimson] (0, 1.5) -- (3.897, 8.250) -- 
	 (6.495, 3.75) -- (.6495, .375)  -- (0, 1.5);
	\draw[thick] (6.495, 3.75) -- (1.299, 3.75);
	 \draw [fill] (0, 1.5) node[anchor=east]{$a$};
	 \draw [fill] (6.495, 3.750) node[anchor=west]{$b$};
	 \draw [fill] (1.299, 3.750) circle [fill,radius=0.1] node[anchor=east]{\hskip-5mm$c_1=c_2$};
	 \draw (3.897, 8.250) node[anchor=west]{$c_3$};
 \end{tikzpicture}
\caption{Polytope GG}
\label{subfig:N=3 GG}
\end{subfigure}}
	\quad
{\begin{subfigure}{0.3\textwidth}
\begin{tikzpicture}[scale=0.38] 
\WeylChamber
	\draw[very thick,fill=Crimson] (0, 3.3) -- (2.858, 8.25) -- 
	 (5.716, 3.3) -- (1.429, .825)  -- (0, 3.3);
	\draw[thick] (0, 3.3) -- (5.716, 3.3);
	 \draw [fill] (0, 3.3) node[anchor=east]{\hskip-5mm$a=c_1$}; 
	 \node at (0,2.6) [anchor=east] {$=c_2$};
	 \node at (5.716, 3.3) [anchor=north]{$b$};
	 \node at (2.858, 8.25) [anchor=west]{$c_3$};
 \end{tikzpicture}
\caption{Polytope FGH}
\label{subfig:N=3 FGH}
\end{subfigure}}
	\quad
{\begin{subfigure}{0.3\textwidth}
\begin{tikzpicture}[scale=0.38] 
\WeylChamber
	\draw[very thick,fill=Crimson] (0, 4.5) -- (2.598, 9) -- 
	 (5.196, 4.5) -- (3.897, 2.25) --  (1.948, 1.125) -- (0, 4.5);
	\draw[thick] (0, 4.5) -- (5.196, 4.5);
	\draw[thick] (1.299, 2.25) -- (3.897, 2.25);	
	 \draw [fill] (0, 4.5) node[anchor=east]{$a=c_1$};
	 \draw [fill] (5.196, 4.5) node[anchor=west]{$b$};
	 \draw [fill] (1.299, 2.25) circle [fill,radius=0.1] node[anchor=east]{$c_2$};
	 \node at  (2.598, 9) [anchor=south]{$c_3$};
 \end{tikzpicture}
\caption{Polytope FH}
\label{subfig:N=3 FH}
\end{subfigure}}

\vskip 5mm

{\begin{subfigure}{0.3\textwidth}
\begin{tikzpicture}[scale=0.38] 
\WeylChamber
	\draw[very thick,fill=Crimson] (0, 6) -- (1.299, 8.25) -- 
	 (3.897, 3.75) --  (2.598, 1.5) -- (0, 6);
	\draw[thick] (3.897, 3.75) -- (1.299, 3.75);
	 \node at (0, 6) [anchor=east]{$a$};
	 \node at  (3.897, 3.75) [anchor=west]{$b$};
	 \draw [fill] (1.299, 3.75) circle [fill,radius=0.1] node[anchor=east]{$c_1$};
	 \node at  (2.598, 1.5) [anchor=north]{$c_2$};
	 \node at  (1.299, 8.25) [anchor=south west]{$c_3$};
 \end{tikzpicture}
\caption{Polytope CH}
\label{subfig:N=3 CH}
\end{subfigure}}
	\quad
{\begin{subfigure}{0.3\textwidth}
\begin{tikzpicture}[scale=0.38] 
\WeylChamber
	\draw[very thick,fill=Crimson] (0, 4.8) -- (2.078, 8.4) -- 
	 (4.157, 4.8) -- (2.078, 1.2)  -- (0, 4.8);
	\draw[thick] (0., 4.8) -- (4.157, 4.8);
	 \node at  (0, 4.8) [anchor=east]{$a=c_1$\!\!};
	 \node at  (4.157, 4.8) [anchor=west]{$b$};
	 \node at  (2.078, 1.2) [anchor=north]{$c_2$};
	 \node at  (2.078, 8.4) [anchor=west]{$c_3$};
 \end{tikzpicture}
\caption{Polytope CEFH}
\label{subfig:N=3 CF}
\end{subfigure}}
	\quad
{\begin{subfigure}{0.3\textwidth}
\begin{tikzpicture}[scale=0.38] 
\WeylChamber
	\draw[very thick,fill=Crimson] (0, 4.2) -- (2.338, 8.25) -- 
	 (4.157, 5.1)  -- (1.819, 1.05) --(0, 4.2);
	\draw[thick] (4.157, 5.1) -- (0.5196, 5.1);
	 \node at  (0, 4.2) [anchor=east]{\hskip5mm$a$};
	 \node at  (4.157, 5.1) [anchor=west]{$b$};
	 \draw [fill] (0.5196, 5.1) circle [fill,radius=0.1] node[anchor=south]{$c_1\;\;$};
	 \node at  (1.819, 1.05) [anchor=north]{$c_2$};
	 \node at  (2.338, 8.25) [anchor=west]{$c_3$};
 \end{tikzpicture}
\caption{Polytope EF}
\label{subfig:N=3 EF}
\end{subfigure}}

\vskip 5mm

{\begin{subfigure}{0.3\textwidth}
\begin{tikzpicture}[scale=0.38] 
\WeylChamber
	\draw[very thick,fill=Crimson] (0, 4.8) -- (2.078, 8.4) -- 
	 (5.196, 3) -- (2.078, 1.2) -- (1.039, 3) -- (0, 4.8);
	\draw[thick] (5.196, 3) -- (1.039, 3);
	 \node at  (0, 4.8) [anchor=east]{$a$};
	 \node at  (5.196, 3.) [anchor=north]{$b$};
	 \draw [fill] (1.039, 3) circle [radius=0.1] node[anchor=east]{\hskip-5mm$c_1=c_2$};
	 \node at  (2.078, 8.4) [anchor=west]{$c_3$};
 \end{tikzpicture}
\caption{Polytope HH}
\label{subfig:N=3 HH}
\end{subfigure}}
	\quad
{\begin{subfigure}{0.3\textwidth}
\begin{tikzpicture}[scale=0.38] 
\WeylChamber
	\draw[very thick,fill=Crimson] (0, 5.1) -- (1.948, 8.475) -- 
	 (3.897, 5.1) -- (1.948, 1.725)  -- (0, 5.1);
	\draw[thick] (0, 5.1) -- (3.897, 5.1);
	 \node at  (0, 5.1) [anchor=east]{\hskip-3mm$a=c_1$\!};
	 \node at  (3.897, 5.1) [anchor=west]{$b$};
	 \node at  (1.948, 1.725) [anchor=east]{$c_2$};
	 \node at  (1.948, 8.475) [anchor=west]{$c_3$};
 \end{tikzpicture}
\caption{Polytope CE}
\label{subfig:N=3 CE}
\end{subfigure}}
	\quad
{\begin{subfigure}{0.3\textwidth}
\begin{tikzpicture}[scale=0.38] 
\WeylChamber
	\draw[very thick,fill=Crimson] (0, 3) -- (3.118, 8.4) -- 
	 (5.716, 3.9)  -- (4.937, 2.85) -- (1.299, .75) -- (0, 3);
	\draw[thick] (0, 3) -- (4.937, 2.85);
	\draw[thick] (5.716, 3.9) -- (.5196, 3.9);
	 \node at  (0, 3) [anchor=east]{$a=c_2$\!\!};
	 \node at  (5.716, 3.9) [anchor=west]{$b$};
	 \draw [fill] (.5196, 3.9) circle [radius=0.1] node[anchor=east]{$c_1\,$};
	 \node at  (3.118, 8.4) [anchor=west]{$c_3$};
 \end{tikzpicture}
\caption{Polytope FG}
\label{subfig:N=3 FG}
\end{subfigure}}

\vskip 5mm
\caption{The remaining transition polytopes---see Fig.~\ref{fig:transition parameters} for notation}
\label{fig:transitions1}
\end{figure}
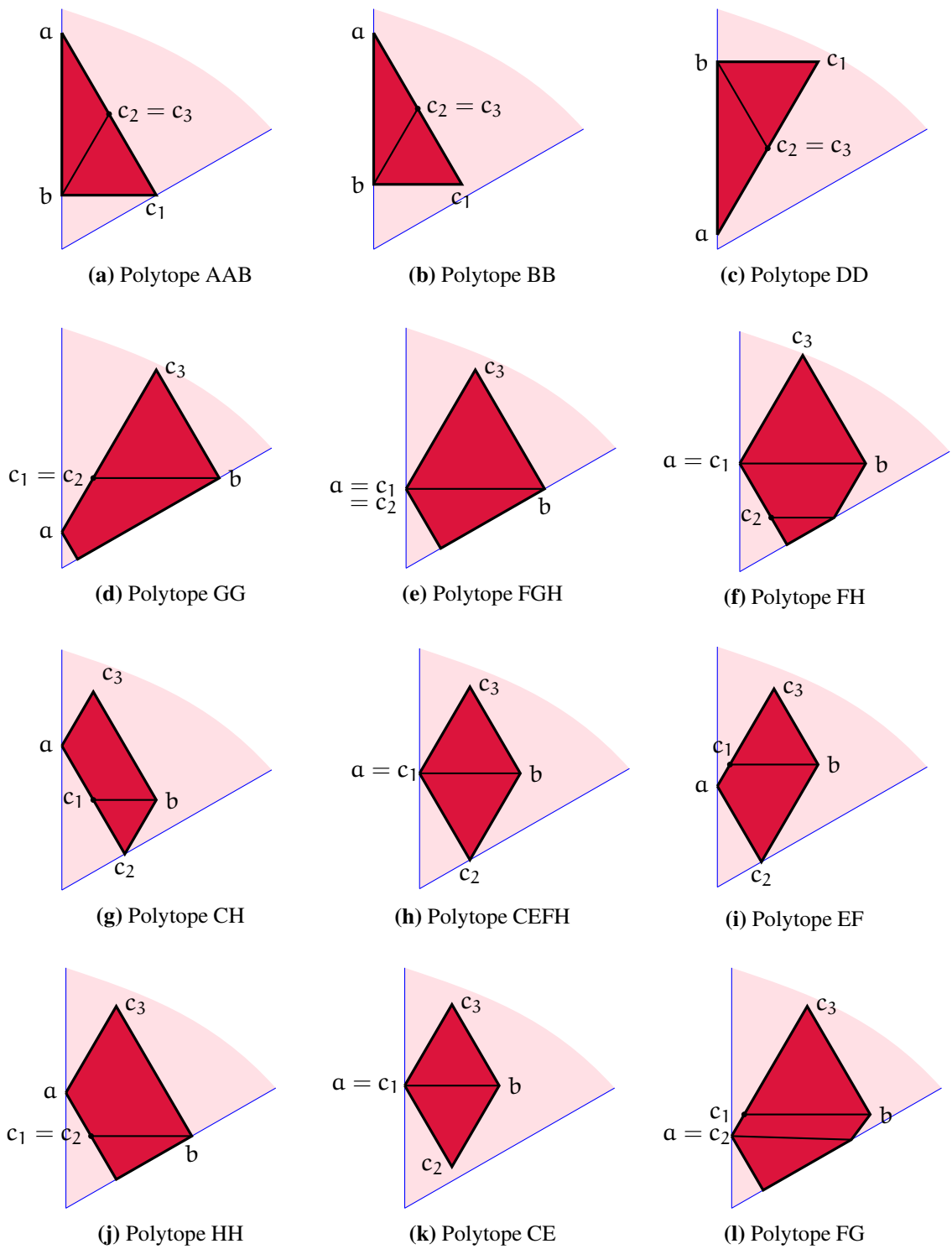

The 8 regions of the $\Gamma$-plane (Figure\,\ref{fig:parameter plane}) are separated by `transition cases', such as occur when $\Gamma_1=\pm\Gamma_2$ or $\Gamma_1=\Gamma_2+\Gamma_3$.  There are also the possibilities of $\Gamma_j=0$, but these we are excluding from consideration as the polytope coincides with that of the 2 point polytopes described in Section\,\ref{sec:polytopes2vortices}. At each of these transitions, one of the vertices of the polytope hits a wall of the Weyl chamber, and a different analysis of the weights is required; indeed the dimension of the symplectic slice is no longer the same.  It can also happen that two of the vertices coincide (such as $c_2=c_3$ when $\Gamma_2=\Gamma_3$ in the transitions denoted BB and DD), but this does not effect  the weight calculations, which are local in $M$.  The transition polytopes are illustrated in Figures\,\ref{fig:transitions2} and \ref{fig:transitions1}; the notation for the different transition cases are shown in Figure\,\ref{fig:transition parameters}.

\begin{figure}
	\centering


{\begin{subfigure}{0.3\textwidth}
\begin{tikzpicture}[scale=0.4] 
\WeylChamber
	\draw[very thick,fill=Crimson] (0,8.25) node[anchor=east]{$a$} -- 
	(1.3,6) circle [fill,radius=0.08] node[black,anchor=south west]{$c_3$} --
	(1.948, 4.875) circle [fill,radius=0.08]  node[black,anchor=south west]{$c_2$} --
	(3.248, 2.625) node[black,anchor=west]{$c_1$} --
	(0.6495, 2.625) node[black,anchor=north]{$b$} --
	(0,3.75) -- (0,8.25);
	\draw [thick] (1.948, 4.875) -- (0.6495, 2.625);
	\draw [thick] (1.3,6) -- (0,3.75);
 \end{tikzpicture}
\caption{Polytope B}
\label{subfig:N=3 B1}
\end{subfigure}}
\quad
{\begin{subfigure}{0.3\textwidth}
\begin{tikzpicture}[scale=0.4] 
\WeylChamber
	\draw[very thick,fill=Crimson] (0, 9) -- (3.897, 2.250) -- 
	 (1.299, 2.25) --  (0, 4.5) -- (0, 9);
	\draw[thick] (2.598, 4.5) -- (1.299, 2.25);
	\draw[thick] (1.299, 6.75) -- (0, 4.5);
	 \draw [fill] (0, 9) node[anchor=east]{$a$};
	 \draw [fill] (1.299, 2.25) node[anchor=north]{$b$};
	 \draw [fill] (3.897, 2.25)  node[anchor=north west]{$c_1$};
	 \draw [fill] (2.598, 4.5) circle [fill,radius=0.1] node[anchor=south west]{$c_2$};
	 \draw [fill] (1.299, 6.750) circle [fill,radius=0.1] node[anchor=south west]{$c_3$};
 \end{tikzpicture}
\caption{Polytope AB}
\label{subfig:N=3 AB}
\end{subfigure}}
\quad
{\begin{subfigure}{0.3\textwidth}
\begin{tikzpicture}[scale=0.4] 
\WeylChamber
	\draw[very thick,fill=Crimson] (0, 8.25) -- (3.572, 2.062) -- 
	 (2.598, 1.5) -- 	(1.299, 1.5) -- (0, 3.75) -- (0, 8.25);
	\draw[thick] (1.299, 6) -- (0, 3.75);
	\draw[thick] (2.598, 3.75) -- (1.3,1.5);
	\draw[thick] (3.248, 2.625) -- (2.598, 1.5);
	 \draw [fill] (0, 8.25) node[anchor=east]{$a$};
	 \draw [fill] (1.3, 1.5) node[anchor=east]{$b$};
	 \draw [fill] (1.299, 6) circle [fill,radius=0.08] node[anchor=south west]{$c_3$};
	 \draw [fill] (2.598, 3.75) circle [fill,radius=0.08] node[anchor=south west]{$c_2$};
	 \draw [fill] (3.248, 2.625) circle [fill,radius=0.08] node[anchor=south west]{$c_1$};
 \end{tikzpicture}
\caption{Polytope A}
\label{subfig:N=3 A1}
\end{subfigure}}

\vskip 5mm
\caption{This shows the transition B $\to$ AB $\to$ A, involving vertex $c_1$ moving to the boundary of the Weyl chamber and getting reflected back but leaving an edge `stuck' to the boundary. See text for further explanation.}
\label{fig:transition A-B}
\end{figure}
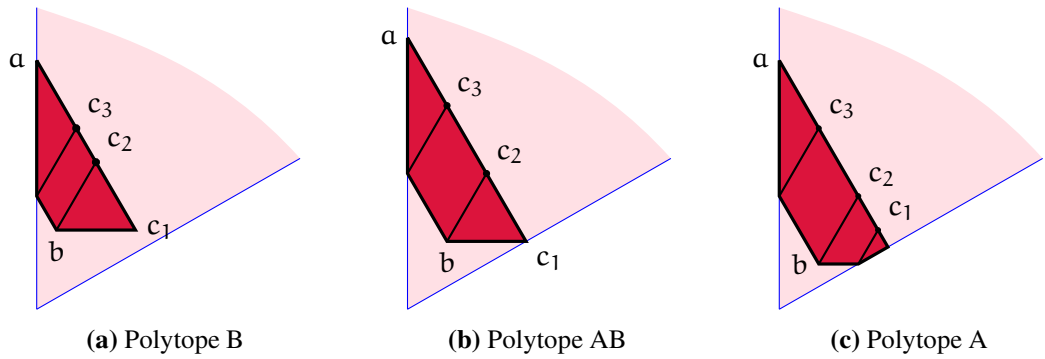

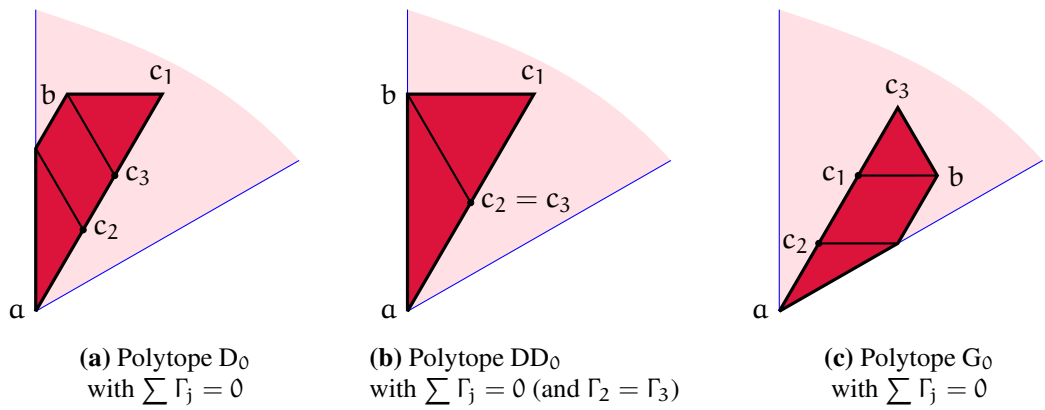
\begin{figure} 
\centering

{\begin{subfigure}{0.3\textwidth}
\begin{tikzpicture}[scale=0.4] 
\WeylChamber
	\draw[very thick,fill=Crimson] (0, 0) -- (4.157, 7.2) -- (1.039, 7.200) --
	 (0, 5.4) -- (0, 0);
	\draw[thick] (1.039, 7.2) -- (2.598, 4.5);
	\draw[thick] (1.559, 2.7) -- (0., 5.4);
	 \draw [fill] (0, 0) node[anchor=east]{$a$};
	 \draw [fill] (1.039, 7.2) node[anchor=east]{$b$};
	 \draw [fill] (4.157, 7.2) node[anchor=south]{$c_1$};
	 \draw [fill] (1.559, 2.7) circle [fill,radius=0.1] node[anchor=west]{$c_2$};
	 \draw [fill] (2.598, 4.5) circle [fill,radius=0.1] node[anchor=west]{$c_3$};
 \end{tikzpicture}
\caption{Polytope D$_0$\\ with $\sum\Gamma_j=0$}
\label{subfig:N=3 D0}
\end{subfigure}}
	\quad
{\begin{subfigure}{0.3\textwidth}
\begin{tikzpicture}[scale=0.4] 
\WeylChamber
	\draw[very thick,fill=Crimson] (0, 0) -- (4.157, 7.2) -- (0., 7.2) --(0,0);
	\draw[thick] (0., 7.2) -- (2.078, 3.6);
	 \draw [fill] (0, 0) node[anchor=east]{$a$};
	 \draw [fill] (0,7.2) node[anchor=east]{$b$};
	 \draw [fill] (4.157, 7.2)  node[anchor=south ]{$c_1$};
	 \draw [fill] (2.078, 3.6) circle [fill,radius=0.1] node[anchor=west]{$c_2=c_3$};
 \end{tikzpicture}
\caption{Polytope DD$_0$\\ with $\sum\Gamma_j=0$ (and $\Gamma_2=\Gamma_3$)}
\label{subfig:N=3 DD0}
\end{subfigure}}
	\quad
{\begin{subfigure}{0.3\textwidth}
\begin{tikzpicture}[scale=0.4] 
\WeylChamber
	\draw[very thick,fill=Crimson] (0, 0) -- (3.897, 6.75) -- 
	 (5.196, 4.5)  -- (3.897, 2.25) --(0, 0);
	\draw[thick] (5.196, 4.5) -- (2.598, 4.5);
	\draw[thick] (1.299, 2.25) -- (3.897, 2.25);
	 \draw [fill] (0, 0) node[anchor=east]{$a$};
	 \draw [fill] (5.196, 4.5) node[anchor=west]{$b$};
	 \draw [fill] (2.598, 4.5) circle [fill,radius=0.1]  node[anchor=east]{$c_1$};
	 \draw [fill] (1.299, 2.25) circle [fill,radius=0.1] node[anchor=east]{$c_2$};
	 \draw [fill] (3.897, 6.750) node[anchor=south]{$c_3$};
 \end{tikzpicture}
\caption{Polytope G$_0$\\ with $\sum\Gamma_j=0$}
\label{subfig:N=3 G0}
\end{subfigure}}

\vskip 5mm
\caption{Polytopes arising for $\Gamma_1+\Gamma_2+\Gamma_3=0$, which implies $a=0$. Notice that D$_0$ and G$_0$ are related by a reflection in the centre line of the Weyl chamber; this is because reversing the signs of the $\Gamma_j$ converts region G$_0$ to D$_0$, via the involution $*$ described in Remark\,\ref{rmk:involution}. A similar observation relates the polytopes for DD$_0$ and GG$_0$ (the latter not drawn). See Figure\,\ref{fig:transition parameters} for the regions in parameter space.}
\label{fig:sum=0}
\end{figure}

However, rather than repeating the weight calculations, it is sufficient to use a continuity argument.  Since $M$ (and hence $M/G$) is compact, and the momentum map depends continuously on $\Gamma$, it follows that the image of the orbit momentum map $\j$ also depends continuously on $\Gamma$. 
It suffices therefore to follow the movement of the vertices as one approaches the boundary of any particular region to conclude the shape of each of the transition momentum polytopes.

In Figure\,\ref{fig:transition A-B} one sees the transition between a polytope of Type A and one of Type B, through the intermediate AB (which occurs when $\Gamma_1=\Gamma_2+\Gamma_3$, assuming as always that $\Gamma_1\geq\Gamma_2\geq\Gamma_3$).  In region B, the point $c_1$ is the image under $J$  of $m=(e_1,e_2,e_2)$. As the $\Gamma_j$ are varied towards region A, the point $J(m)$ moves towards the boundary of the positive Weyl chamber, and crosses the wall
so that in region A the point $J(m)$ is no longer in the positive Weyl chamber. In region A, the point $c_1$ is the value $J(m')$ for $m'=(e_2,e_1,e_1)$. At the transition (Type AB), both of these points map to the same point in the wall of $\tt^*_+$.  Note that in general, the two values $J(m)$ and $J(m')$ have the same spectrum and are related by an element of the Weyl group.

\paragraph{Symplectic reduction}
In the companion paper \cite{MS19b}, we describe the reduced spaces $M_\mu$ for $\mu\in\Delta(M)$.  The possibilities for $M_\mu$ are a sphere, a sphere with singularities and a single point.

\small
\setlength{\parskip}{0pt}

\bigskip
\setlength{\parindent}{0pt}

JM: j.montaldi@manchester.ac.uk

AS: amna.shaddad@gmail.com

\medskip

School of Mathematics\\
University of Manchester\\
Manchester M13 9PL, UK

\end{document}